\documentclass[12 pt]{amsart}
\usepackage{a4}
\usepackage[french]{babel}
\usepackage[latin1]{inputenc}
\usepackage{xypic}
\xyoption{all}
\usepackage{epsfig}

\newtheorem{thm}{Th\'eor\`eme}[section]
 
\newtheorem{lem}[thm]{Lemme}
  
\newtheorem{cor}[thm]{Corollaire}
\newtheorem{prop}[thm]{Proposition}

\theoremstyle{definition}
\newtheorem{defn}[thm]{D\'efinition}
\newtheorem{exmp}[thm]{Exemple}

\theoremstyle{remark}
\newtheorem{rem}[thm]{Remarque} 

\newcommand {\spec}{{\rm Spec }}

\newcommand {\ord}{{\rm ord}}

\newcommand {\card}{{\rm card}}

\newcommand {\Div}{{\rm div}}
\newcommand {\som}{{\rm som}}
\newcommand {\aret}{{\rm ar}}

\newcommand {\pic}{{\rm Pic}}
\newcommand {\kuro}{\displaystyle{\frac{R[[Z_1,Z_2]]}{(Z_1Z_2-\pi^e)}}}

\begin{document}

\title[Arbres de Hurwitz]
{Arbres de Hurwitz et automorphismes d'ordre $p$ des disques et 
des couronnes $p$-adiques formels}
\author{Yannick Henrio}
\date{28 janvier 2000}
\maketitle
\noindent
\begin{center}
Laboratoire de Math\'ematiques Pures de Bordeaux\\
UPRES-A 5467 CNRS\\
Universit\'e Bordeaux I\\
351 cours de la lib\'eration 33 405 Talence cedex, France\\
fax : (33) 05 56 84 69 29\\
{\texttt e-mail : henrio@math.u-bordeaux.fr}
\end{center}\ \\
\newline\newline
\par
Soit $R$ un anneau de valuation discr\`ete complet, de corps des 
fractions $K$ de caract\'eristique nulle et de corps r\'esiduel 
alg\'ebriquement clos $k$, de caract\'eristique $p>0$. On notera 
$\pi$ une uniformisante de $R$ et $v_K$ la valuation de $K$, normalis\'ee par 
$v_K(\pi)=1$. 
On supposera en outre que $K$ contient une racine primitive $p$-i\`eme 
$\zeta=\lambda+1$ 
de l'unit\'e. Le disque formel sur $R$ est le $R$-sch\'ema 
$\mathcal D:=\spec R[[Z]]$. Si $e$ est un entier strictement positif sur $R$, 
la couronne formelle d'\'epaisseur $e$ est le $R$-sch\'ema 
$\mathcal C_e:=\spec \kuro$. De tels $R$-sch\'emas se rencontrent dans 
l'\'etude des courbes relatives sur $R$ : 
En effet, soit $X$ un $R$-sch\'ema plat, dont les 
fibres sont de dimension $1$. Si la fibre g\'en\'erique de $X$ 
est lisse sur $K$, pour un point ferm\'e $x$ de la fibre sp\'eciale $X_k$ 
de $X$, le $R$-sch\'ema local $\spec \hat {\mathcal O_{X,x}}$ est un disque 
formel si $x$ est un point lisse de $X_k$ et une couronne formelle si $x$ 
un point double ordinaire. Les probl\`emes de rel\`evements galoisiens 
(voir \cite{G-M 1} et \cite{He1}) 
pour les courbes conduisent alors \`a des questions d'existence 
d'automorphismes de disques ou de couronnes formels.
\par
Dans \cite{G-M 2}, B. Green et M. Matignon ont \'etudi\'e 
les automorphismes d'ordre $p$ du disque formel sur $R$. Ils montrent 
notamment que les positions relatives des points fixes de la fibre 
g\'en\'erique sous l'action d'un tel automorphisme $\sigma$ satisfont des 
contraintes fortes, de nature m\'etrique et alg\'ebrique. Plus 
pr\'ecis\'ement, ils montrent tout d'abord l'existence et l'unicit\'e 
de l'\'eclatement minimal de $\mathcal D$, \`a support dans la fibre 
sp\'eciale, pour lequel les sp\'ecialisations des points fixes de 
l'automorphisme $\sigma$ sont lisses et distinctes. 
Sa fibre sp\'eciale est un arbre de 
droites projectives sur $k$, orient\'e \`a partir de la transform\'ee 
stricte de $\mathcal D \times_R k$. Ils montrent alors que les points fixes se 
sp\'ecialisent dans les bouts de l'arbre ; de plus, sur un tel bout $E$, 
il existe une unique forme diff\'erentielle logarithmique dont le diviseur 
a pour unique z\'ero le point double et pour uniques p\^oles les 
sp\'ecialisations des points fixes 
sur $E$, et les r\'esidus en ces points sont reli\'es directement \`a l'action 
sur l'espace tangent aux points fixes. Par ailleurs, sur une composante 
interne de l'arbre, il existe une forme diff\'erentielle exacte 
dont le diviseur a pour support les points doubles de cette composante.
\par
Dans cet article, nous introduisons la notion combinatoire 
d'arbre de Hurwitz, qui permet 
d'exprimer ces contraintes de fa\c con synth\'etique. Cette notion est 
adapt\'ee \'egalement au cas des automorphismes d'ordre $p$ de 
couronnes formelles, et nous \'etendons les r\'esultats de Green et 
Matignon dans ce contexte. Nous donnons par ailleurs une condition 
n\'ecessaire et suffisante pour qu'un arbre de Hurwitz provienne d'un 
automorphisme d'ordre $p$ d'un disque ou d'une couronne formelle. Dans le 
cas du disque, on obtient ainsi la r\'eciproque du th\'eor\`eme III 2.1 
de \cite{G-M 2}.
\par
Dans une premi\`ere partie, nous rappellons des r\'esultats connus 
sur la r\'eduction des torseurs de Kummer de degr\'e $p$. Le seul r\'esultat 
essentiellement nouveau est le corollaire \ref{reducbord}, qui permet 
d'\'ecrire explicitement l'action d'un automorphisme d'ordre $p$ au bord 
d'un disque formel. La seconde partie est consacr\'ee \`a la d\'efinition 
d'un arbre de Hurwitz. La troisi\`eme partie traite des automorphismes 
d'ordre $p$ au bord d'un disque formel ; nous construisons l'arbre de 
Hurwitz associ\'e et, \`a l'aide de la g\'eom\'etrie formelle, nous montrons 
le {\bf Th\'eor\`eme de r\'ealisation :}\\
\newline
{\it \par Soit $(\Gamma, \mathcal H)$ un arbre de Hurwitz, il 
provient d'un $R$-automorphisme d'ordre $p$ de $R[[Z]]$ si et seulement si 
il v\'erifie les conditions $D[i]$ ci-dessous ($1 \le i \le 3$) :
\newline
$D[1]$ La racine de $(\Gamma, \mathcal H)$ est de valence $1$.
\newline
$D[2]$ Toute ar\^ete aboutissant \`a un sommet maximal est une feuille.
\newline
$D[3]$ Tout sommet de valence sup\'erieure ou \'egale \`a $3$ est 
r\'ealisable.}
\newline\newline
Nous invitons le lecteur \`a entrer plus avant dans la lecture 
pour comprendre la terminologie employ\'ee. La condition $D[3]$ consiste en 
l'existence de certaines formes diff\'erentielles. Nous compl\'etons ensuite 
ce r\'esultat par un crit\`ere suffisant, mais non n\'ecessaire, pour assurer 
l'existence de ces formes diff\'erentielles. Enfin, dans la derni\`ere partie, 
nous donnons des r\'esultats similaires pour les couronnes formelles, 
ainsi qu'un th\'eor\`eme de structure pour les arbres de Hurwitz associ\'es 
\`a un automorphisme d'ordre $p$ \`a petits conducteurs aux bords, analogue 
au th\'eor\`eme III 3.1 de \cite{G-M 2}.
\par
Nous noterons $R\{T_1,\dots,T_n\}$ l'alg\`ebre des s\'eries restreintes \`a 
coefficients dans $R$, c'est-\`a-dire l'alg\`ebre des s\'eries formelles 
$$f=\sum_{(\nu_1,\dots,\nu_n) \in \mathbb N^n}
a_{\nu_1,\dots ,\nu_n} T_1^{\nu_1}\dots T_n^{\nu_n}$$
avec $\displaystyle{\lim_{\nu_1+\dots +\nu_n \to +\infty}
a_{\nu_1,\dots ,\nu_n}=0}$.
\par
Nous noterons $R[[T]]\{T^{-1}\}$ l'alg\`ebre des s\'eries de Laurent 
$f=\sum_{\nu \in \mathbb Z}a_{\nu}T^{\nu}$ telles que 
$\displaystyle{\lim_{\nu \to -\infty}a_{\nu}=0}$. C'est un anneau de 
valuation discr\`ete complet d'uniformisante $\pi$ et de corps 
r\'esiduel $k((T))$. 
Si $\mathcal X$ est un disque ou une couronne formelle, on appelle bord de 
$\mathcal X$ un point g\'en\'erique d'une composante irr\'eductible de la 
fibre sp\'eciale (qui est un diviseur de Weil sur $\mathcal X$). 
Le sch\'ema $\mathcal X$ \'etant normal, un bord $\eta$ d\'efinit 
une valuation not\'ee $v_{\eta}$ du corps des fonctions de $\mathcal X$. 
On v\'erifie que pour un bord $\eta$ de $\mathcal X$, 
on a un isomorphisme de $R$-alg\`ebre de $\hat {\mathcal O_{\mathcal X,\eta}}$ 
sur $R[[T]]\{T^{-1}\}$ ; un \'el\'ement $T$ de $\hat {\mathcal O_{\mathcal X,\eta}}$ tel que $\hat{\mathcal O_{\mathcal X,\eta}}=R[[T]]\{T^{-1}\}$ sera 
appel\'e param\`etre au bord $\eta$. Le corps r\'esiduel $k((T))$ de 
$\hat {\mathcal O_{\mathcal X,\eta}}$ est muni 
d'une valuation discr\`ete $\ord_{\eta}$, normalis\'ee de fa\c con \`a ce 
qu'une uniformisante de $k((T))$ ait pour valuation $1$. 
Si $f$ est non nul dans $\hat {\mathcal O_{\mathcal X,\eta}}$, 
$f$ s'\'ecrit de mani\`ere unique $f=\pi^{v_{\eta}(f)}f_0$, avec 
$f_0$ non nul. On note alors 
$\ord_{\eta}(f):=\ord_{\eta}(\bar f_0)$, o\`u $\bar f_0$ est 
l'image r\'esiduelle de $f_0$ dans $k((T))$.

\section{R\'eduction des torseurs de Kummer de degr\'e $p$}
\subsection{Les sch\'emas en groupes $\mathcal G^{(n)}$ et 
$\mathcal H_n$}
\par
On rappelle ici la d\'efinition de certains sch\'emas en groupes qui 
interviennent dans la r\'eduction des $\mu_p$-torseurs. Pour plus 
d'informations concernant ces sch\'emas en groupes et leurs applications \`a 
la d\'eformation d'Artin-Schreier \`a Kummer, on renvoie \`a \cite{OSS}. 
Pour tout entier $n$ strictement positif, on note 
$\mathcal G^{(n)}:=\displaystyle{\spec R[x,\frac{1}{\pi^nx+1}]}$ le 
sch\'ema en groupes affine sur $R$ dont l'alg\`ebre de Hopf est donn\'ee 
ci-dessous : 
\begin{itemize}
\item
comultiplication : 
\begin{eqnarray*}
R[x,\frac{1}{\pi^nx+1}] & & \to R[x,\frac{1}{\pi^nx+1}] 
\otimes_R R[x,\frac{1}{\pi^nx+1}]\\
 x & & \mapsto  x\otimes 1 + 1 \otimes x +\pi^n x \otimes x
\end{eqnarray*}
\item
coinverse :
\begin{eqnarray*}
R[x,\frac{1}{\pi^nx+1}] & & \to R[x,\frac{1}{\pi^nx+1}]\\
x & & \mapsto  -\frac x{\pi^nx+1}
\end{eqnarray*}
\item
coidentit\'e
\begin{eqnarray*}
R[x,\frac{1}{\pi^nx+1}] & & \to R\\
x & & \mapsto  0
\end{eqnarray*}
\end{itemize}
\par
La fibre sp\'eciale de $\mathcal G^{(n)}$ s'identifie canoniquement au groupe 
additif $\mathbb G_a$ sur $k$ ; par aillleurs, la fibre g\'en\'erique 
est isomorphe au groupe multiplicatif $\mathbb G_m$ sur $K$ : Plus 
pr\'ecis\'ement, on d\'efinit un homomorphisme de sch\'emas en groupes 
$\phi_n$ de 
$\mathcal G^{(n)}$ dans $\mathbb G_{m,R}$ par 
l'homomorphisme de $R$-alg\`ebres de Hopf :
\begin{eqnarray*}
R[u,u^{-1}] & & \to R[x,\frac{1}{\pi^nx+1}]\\
u & & \mapsto \pi^nx+1
\end{eqnarray*}
On a alors la suite exacte pour la 
topologie $fppf$ :
$$(1) \quad 0 \to \mathcal G^{(n)} \stackrel{\phi_n}{\longrightarrow} 
\mathbb G_{m,R} \to j_{n*}\mathbb G_{m,R_n} \to 0$$
o\`u $R_n:=R/(\pi^n)$ et $j_n$ est l'immersion ferm\'ee de $\spec R_n$ dans 
$\spec R$. (voir \cite{OSS}). On voit alors que 
$\phi_n \otimes_R K$ est un isomorphisme.
\par
Pour $0<n \le v_K(\lambda)=\displaystyle{\frac{v_K(p)}{p-1}}$, le polyn\^ome 
$\displaystyle{\frac{(\pi^nx+1)^p-1}{\pi^{pn}}}$ est \`a coefficients dans 
$R$ ; on d\'efinit alors un 
homomorphisme de $R$-sch\'emas en groupes $\psi_n$ de 
$\mathcal G^{(n)}$ dans $\mathcal G^{(pn)}$ par :
\begin{eqnarray*}
R[y,\frac{1}{\pi^{pn}y+1}] & & \to R[x,\frac{1}{\pi^nx+1}]\\
y & & \mapsto  \frac{(\pi^nx+1)^p-1}{\pi^{pn}}
\end{eqnarray*}
\par
L'homomorphisme $\psi_n$ est une isog\'enie de degr\'e $p$, on note 
$\mathcal H_n$ le noyau de $\psi_n$, et on obtient la suite exacte pour la 
topologie $fppf$ :
$$(2) \quad 0 \to \mathcal H_n \to \mathcal G^{(n)} 
\stackrel{\psi_n}{\longrightarrow} 
\mathcal G^{(pn)} \to 0$$
Le sch\'ema en groupes $\mathcal H_n$ est fini et plat sur $R$, 
de degr\'e $p$. Sa fibre g\'en\'erique est isomorphe au groupe $\mu_{p,K}$. 
Si $0<n<v_K(\lambda)$, sa fibre sp\'eciale est le groupe radiciel additif 
$\alpha_p$ ; si $n=v_K(\lambda)$, sa fibre sp\'eciale est un groupe \'etale 
isomorphe \`a $\mathbb Z/p\mathbb Z$.
\subsection{Torseurs $fppf$ sous $\mathcal H_n$}
\par
Soit $n$ un entier strictement positif, $X$ un $R$-sch\'ema, on d\'eduit de 
la suite exacte $(1)$ la suite exacte 
$$H^0(X,\mathbb G_m) \to H^0(X_n,\mathbb G_m) \to H^1(X,\mathcal G^{(n)}) \to 
\pic X\to \pic X_n$$
o\`u $X_n:=X \times_{\spec R}\spec R_n$. En particulier, on a le 
\begin{lem}
Si $\pic X \to \pic X_n$ est injectif et si 
$H^0(X,\mathbb G_m) \to H^0(X_n,\mathbb G_m)$ est 
surjectif, $H^1(X,\mathcal G^{(n)})=0$. En particulier, si $X=\spec A$, avec 
$A$ une $R$-alg\`ebre compl\`ete pour la topologie 
$\pi$-adique, $H^1(X,\mathcal G^{(n)})=0$.
\end{lem}
\par
Soit $X$ un $R$-sch\'ema, si $0<n\le v_K(\lambda)$, on d\'eduit de la suite 
exacte $(2)$ la suite exacte 
$$H^0(X,\mathcal G^{(n)}) \to H^0(X,\mathcal G^{(pn)}) \to H^1(X,\mathcal H_n) \to 
H^1(X,\mathcal G^{(n)})$$
\begin{cor}
Si $X=\spec A$, avec 
$A$ une $R$-alg\`ebre compl\`ete pour la topologie $\pi$-adique,
$$H^1(X,\mathcal H_n)=\mathcal G^{(pn)}(A)/\mathcal G^{(n)}(A).$$
\end{cor}
\begin{rem} On peut concr\`etement re\'ecrire le corollaire de la mani\`ere 
suivante : On consid\`ere une $R$-alg\`ebre $A$ factorielle et compl\`ete 
pour la topologie $\pi$-adique, un torseur sous $\mathcal H_n$ au-dessus de 
$\spec A$ est alors 
de la forme $\spec B$, avec
$$B:=\displaystyle{\frac{A[w]}{\frac{(\pi^nw+1)^p-1}{\pi^{pn}}-u}},$$
o\`u $u$ appartient \`a $A$. Deux \'el\'ements $u_1$ et 
$u_2$ de $A$ d\'efinissent le m\^eme torseur s'il existe $v$ dans $A$ tel que 
$$u_2=u_1(\pi^nv+1)^p+\frac{(\pi^nv+1)^p-1}{\pi^{pn}}.$$
\end{rem}
\subsection{Torseurs $fppf$ sous $\mu_p$}
\par
On a la suite exacte de sch\'emas en groupes commutatifs, pour la topologie 
$fppf$ :
$$1 \to \mu_p \to \mathbb G_m \stackrel{x \mapsto x^p}
{\longrightarrow} \mathbb G_m \to 1.$$
On en d\'eduit le
\begin{lem}
Soit $X$ un $R$-sch\'ema, si $\pic X=0$, 
$H^1(X,\mu_p)=H^0(X,\mathbb G_m)/H^0(X,\mathbb G_m)^p$. En particulier, 
si $X=\spec A$, avec 
$A$ une $R$-alg\`ebre factorielle,
$$H^1(X,\mu_p)=\mathbb G_m(A)/\mathbb G_m(A)^p.$$
\end{lem}
\begin{rem} Comme ci-dessus, on peut reformuler concr\`etement le lemme 
de la mani\`ere suivante : On consid\`ere une $R$-alg\`ebre $A$ factorielle, 
un torseur sous $\mu_p$ au-dessus de 
$\spec A$ est alors 
de la forme $\spec B$, avec
$$B=\frac{A[y]}{(y^p-u)}.$$
o\`u $u$ est une unit\'e de $A$, unique \`a la multiplication pr\`es 
d'une puissance $p$-i\`eme d'une unit\'e de $A$.
\end{rem}

\subsection{R\'eduction des $\mu_p$-torseurs}
\begin{prop}
\label{reducmuptorseur}
\par
Soit $X:=\spec A$ un sch\'ema affine plat sur $R$, dont les fibres sont 
int\`egres et de dimension $1$ ; on suppose que $A$ est 
une $R$-alg\`ebre factorielle et compl\`ete pour la topologie 
$\pi$-adique. 
Soit $Y_K \to X_K$ un $\mu_p$-torseur \'etale non trivial, donn\'e par 
une \'equation $y^p=f$, o\`u $f$ est 
inversible dans $A_K$, et $Y$ le normalis\'e de 
$X$ dans $Y_K$ ; on suppose que 
la fibre sp\'eciale de $Y$ est int\`egre. Soit $\eta$ (resp. $\eta'$) 
le point g\'en\'erique de la fibre sp\'eciale de $X$ (resp. $Y$). Les anneaux 
locaux $\mathcal O_{X,\eta}$ et $\mathcal O_{Y,\eta'}$ sont alors des anneaux 
de valuation discr\`ete d'uniformisante $\pi$. 
Notons $\delta$ la valuation de la 
diff\'erente de $\mathcal O_{Y,\eta'}/\mathcal O_{X,\eta}$. On distingue 
alors deux cas suivant la valeur de $\delta$.
\begin{itemize}
\item
Si $\delta=v_K(p)$, $Y$ est un $\mu_{p,R}$-torseur pour la topologie $fppf$, 
donc $\mathcal Y=\spec B$, avec 
$B:=\displaystyle{\frac{A[y]}{(y^p-u)}}$, o\`u $u$ est une unit\'e de 
$A$, unique \`a la multiplication d'une puissance $p$-i\`eme d'une unit\'e de 
$A$ pr\`es. On dit que le torseur $Y \to X$ a {\bf r\'eduction multiplicative}.
\item
Si $0 \le \delta<v_K(p)$, on a $\delta=v_K(p)-n(p-1)$, o\`u $n$ est un entier 
tel que $0<n\le v_K(\lambda)$, et $Y \to X$ est un torseur sous 
$\mathcal H_n$  pour la topologie $fppf$, donc est donn\'e par $Y=\spec B$, 
$$B:=\displaystyle{\frac{A[w]}
{(\frac{(\pi^nw+1)^p-1}{\pi^{pn}}-u)}},$$
o\`u $u$ est un \'el\'ement de $A$. 
De plus, si $B$ est isomorphe \`a $\displaystyle{\frac{A[w]}
{\frac{(\pi^nw+1)^p-1}{\pi^{pn}}-u'}}$, alors il existe $v$ dans 
$A$ tel que 
$u'=u(\pi^nv+1)^p+\displaystyle{\frac{(\pi^nv+1)^p-1}{\pi^{pn}}}.$\\
Si $0<\delta<v_K(p)$ (resp. $\delta=0$), on dit que le torseur $Y \to X$ a 
{\bf r\'eduction additive} (resp. {\bf r\'eduction \'etale}).
\end{itemize}
\end{prop}
\begin{proof}
\par
Notons $Y:=\spec B$, o\`u $B$ est une 
$A$-alg\`ebre finie, normale. 
On peut \'ecrire $f=\pi^rf_0$, o\`u $r$ est un entier, 
$f_0$ appartient \`a $A$ et n'est 
pas divisible par 
$\pi$ ; on peut en outre supposer $0 \le r <p$. 
Supposons $r>0$, alors comme $y$ appartient \`a $B$, et 
$B \otimes_R k$ est r\'eduit, 
$y=\pi^sy'$, avec $y'$ dans $B$ et $s>0$. Mais alors, $f_0=\pi^{sp-r}{y'}^p$, 
donc $f_0$ appartient \`a $\pi B$. Comme $\displaystyle{\frac {f_0}{\pi}}$ 
est entier sur $A$, en consid\'erant un polyn\^ome unitaire de degr\'e $d$ 
qui l'annulle, on voit que $\bar f_0^d$ est nul. Comme $A\otimes_R k$ 
est int\`egre, l'image $\bar f_0$ de $f_0$ dans $A\otimes_R k$ est nulle, 
i.e. $f_0$ appartient \`a $\pi A$, ce qui contredit sa 
d\'efinition. 
Ainsi, $r=0$ et donc $f$ appartient \`a $A$, $\bar f\ne 0$. En fait, comme 
$\pi$ est irr\'eductible dans $A$ et $f$ est inversible dans $A_K$, $f$ est 
inversible dans $A$.
\par
On suppose tout d'abord que $\bar f$ n'est pas une puissance $p$-i\`eme, alors 
$\displaystyle{\frac{A[y]}{(y^p-f)}}$ est int\'egralement clos, contenu dans 
$B$  et de m\^eme corps des fractions, ainsi il est \'egal \`a 
$B$. On a alors 
$\delta=v_K(p)$ et $\spec B$ est un torseur sous $\mu_p$.
\par
On suppose \`a pr\'esent que $f$ est une puissance $p$-i\`eme modulo $\pi$. Il 
existe donc un $g$ dans $A$, n\'ecessairement inversible, tel que 
$f=g^p \mod \pi$. Rempla\c cant $f$ par $g^{-p}f$, on peut alors supposer 
que $f=1+\pi^{r_0}h_0$, o\`u $h_0$ appartient \`a $A$ et $\bar h_0 \ne 0$. 
Si on avait $r_0>pv_K(\lambda)$, alors $1+\pi^{r_0}h_0$ serait une puissance 
$p$-i\`eme dans $A$ et le torseur $Y_K \to X_K$ serait trivial. Ainsi 
$r_0 \le pv_K(\lambda)$. Si $r_0=pv_K(\lambda)$, $Y$ est un torseur (\'etale) 
sous $\mathcal H_{v_K(\lambda)}$. Sinon, on \'ecrit la division euclidienne 
$r_o=pq_0+s_0$, $0 \le s_0<p$. En posant $y:=\pi^{q_0}y_0+1$, on obtient 
$$\frac{(\pi^{q_0}y_0+1)^p-1}{\pi^{pq_0}}=\pi^{s_0}h_0.$$
Supposons $s_0>0$, on a alors $y_0=\pi y_0'$, avec $y_0'$ dans $B$. De la 
relation\\
$1+\pi^{s_0}h_0=(1+\pi^{q_0+1}y_0')^p$, on tire alors la 
contradiction $\bar h_0=0$. Ainsi, $r_0=pq_0$, avec $q_0<v_K(\lambda)$. Si 
$\bar h_0$ n'est pas une puissance $p$-i\`eme, on a fini. Sinon, on peut 
r\'eappliquer le m\^eme processus et on tombe apr\`es un nombre fini 
d'\'etapes sur une \'equation $y^p=1+\pi^{pn}h$, avec ou bien 
$n=v_K(\lambda)$ et on a une r\'eduction \'etale,  ou bien 
$n<v_K(\lambda)$, $\bar h$ n'est pas une puissance $p$-i\`eme et on a un 
torseur sous $\mathcal H_n$, donc radiciel en r\'eduction.
\end{proof}
\begin{exmp}
Soit $A:=R\{T,(T-T_i)^{-1}\}_{1 \le i\le r}$, $X:=\spec A$ v\'erifie les 
conditions de la proposition pr\'ec\'edente. Un $\mu_p$-torseur au-dessus 
de sa fibre g\'en\'erique se prolonge donc en un torseur sous $\mu_p$ 
(on dira qu'on a r\'eduction multiplicative) ou sous $\mathcal H_n$ 
(on dira qu'on a r\'eduction additive si $n<v_K(\lambda)$ et \'etale si 
$n=v_K(\lambda)$).
\end{exmp}
\par
Le corollaire suivant est essentiel pour la construction d'automorphismes 
d'ordre $p$ du disque formel \`a partir d'un arbre de Hurwitz. Il permet en 
effet d'\'ecrire explicitement, sur un param\`etre convenable, l'action d'un 
tel automorphisme au bord du disque.
\begin{cor}\label{reducbord}
\par
Soit $A:=R[[T]]\{T^{-1}\}$, $B$ une 
$A$-alg\`ebre finie de degr\'e $p$, plate sur $R$ et $G$ un groupe 
$p$-cyclique d'automorphismes de $B$ ; on 
suppose que $A$ s'identifie \`a $B^G$, et que $\bar B:=B\otimes_R k$ 
est un corps, 
extension purement ins\'eparable de degr\'e $p$ de $\bar A=k((t))$. 
L'anneau $B$ est alors un 
anneau de valuation discr\`ete complet d'uniformisante $\pi$, de corps 
r\'esiduel de la forme $k((z))$ ; si $Z$ rel\`eve $z$ dans $B$, 
$B=R[[Z]]\{Z^{-1}\}$. Soit $\delta$ la valuation de 
la diff\'erente de $B$ sur $A$.
\newline
(A) Si $\delta=v_K(p)$, $\spec B\to \spec A$ est un $fppf$-torseur sous 
$\mu_p$, et ainsi il existe une unit\'e $u$ dans $A$ tel que 
$B=\displaystyle{\frac{A[y]}{(y^p-u)}}$. 
De plus, $u$ est unique \`a la multiplication d'une puissance $p$-i\`eme 
pr\`es d'une unit\'e de $A$. En particulier, la forme diff\'erentielle 
$\omega:=\displaystyle{\frac{d\bar u}{\bar u}}$ est uniquement d\'etermin\'ee. 
De plus, le r\'esidu $h$ de $\omega$ appartient \`a $\mathbb F_p$. 
\par
Notons $\phi$ l'unique isomorphisme de $G$ sur 
$\mathbb Z/p\mathbb Z$ tel que pour $\sigma \in G$, on ait 
$\sigma y=\zeta^{\phi(\sigma)}y$, on a alors deux cas :
\newline
a) Si le r\'esidu $h$ de $\omega$ est non nul, alors $y=Z^h$, o\`u $Z$ est un 
param\`etre de $B$, et, pour $\sigma \in G$, 
$\sigma Z=\zeta^{\frac 1h \phi(\sigma)}Z$.
\newline
b) Si le r\'esidu $h$ de $\omega$ est nul, alors l'entier 
$m:=1+\ord_t \omega$ est positif et premier \`a $p$. De plus, (pour 
un choix convenable de $u$) il existe 
un param\`etre $Z$ de $B$ tel que $y=1+Z^m$ et, pour $\sigma \in G$,
$$\sigma Z=\zeta^{\frac 1m \phi(\sigma)}Z
(1+\zeta^{-\phi(\sigma)}(\zeta^{\phi(\sigma)}-1)Z^{-m})
^{\frac 1m}$$
(B) 
Si $0<\delta<v_K(p)$, $\spec B\to \spec A$ est un $fppf$-torseur sous 
$\mathcal H_n$, 
o\`u $n$ est un entier donn\'e par $v_K(p)=\delta+n(p-1)$. Ainsi, 
il existe $u$ dans $A$ tel que 
$B=\displaystyle{\frac{A[w]}{(\frac{(\pi^nw+1)^p-1}{\pi^{pn}}-u)}}$, 
$\bar u \notin k((t^p))$. De plus, $\bar u$ est unique \`a l'addition 
d'une puissance $p$-i\`eme pr\`es. En particulier, la forme diff\'erentielle 
$\omega:=d\bar u$ est uniquement d\'etermin\'ee. L'entier 
$m:=1+\ord_t(\omega)$ est premier \`a $p$. 
\par
Notons $\phi$ l'unique 
isomorphisme de $G$ sur $\mathbb Z/p\mathbb Z$ tel que pour $\sigma \in G$, 
on ait 
$\sigma(\pi^nw+1)=\zeta^{\phi(\sigma)}(\pi^nw+1)$. On a alors $w=Z^m$, o\`u 
$Z$ est un param\`etre de $B$, et pour $\sigma \in G$, 
$\sigma Z=\zeta^{\frac 1m \phi(\sigma)}Z(1+\pi^{-n}(\zeta^{\phi(\sigma)}-1)
\zeta^{-\phi(\sigma)}Z^{-m})^{\frac 1m}$.
\end{cor}
\begin{proof}
Comme $k$ est alg\'ebriquement clos, $\bar B$ est un corps $k((z))$ 
de s\'eries formelles en une variable sur $k$. Soit $Z$ dans $B$ relevant $Z$, 
il suit du lemme 2.1 de \cite{He2} que $B=R[[Z]]\{Z^{-1}\}$. On est donc 
dans la situation de la proposition \ref{reducmuptorseur} (Remarquer que 
$A$ est principal, donc factoriel).
\par
(A) On suppose tout d'abord que $\delta=v_K(p)$. La proposition 
\ref{reducmuptorseur} montre alors que $\spec B\to \spec A$ est un 
$fppf$-torseur sous $\mu_p$, donn\'e par 
$$B:=\frac{A[y]}{(y^p-u)}.$$
(a) Si $\ord_t \bar u$ est premier \`a $p$, alors il est \'egal modulo $p$ 
au r\'esidu $h$ de $\omega:=\displaystyle{\frac{d\bar u}{\bar u}}$. En 
particulier, ce r\'esidu est un inversible de $\mathbb F_p$. De plus, 
$y$ est modulo $\pi$ la puissance $h$-i\`eme d'une uniformisante $z$ 
de $\bar B$. 
Comme $B$ est hens\'elien, $y=Z^h$, o\`u $Z$ rel\`eve $z$ dans $B$, ainsi 
$Z$ est un param\`etre de $B$. Le calcul de $\sigma Z$ est imm\'ediat.\\
(b) Si maintenant $\ord_t \bar u$ est divisible par $p$, on peut le 
supposer nul. Alors, si $\omega:=\displaystyle{\frac{d\bar u}{\bar u}}$, on a 
$1+\ord \omega=1+\ord d\bar u$ positif et premier \`a $p$. En particulier, 
le r\'esidu de $\omega$ est nul. Alors on peut \'ecrire 
$\bar u=1+\sum_{1\le i \le [\frac mp]}\bar c_i^pt^{ip}+\bar at^m \mod t^{m+1}$, o\`u les $c_i$ et $a$ sont dans $R$. 
En multipliant $Y$ par $1-c_1T$, on se ram\`ene \`a $\bar c_1=0$. 
En appliquant 
successivement ce processus, on peut supposer que $\bar u$ s'\'ecrit 
$\bar u=1+\bar at^m \mod t^{m+1}$, avec $\bar a\ne 0$. On a 
$(\bar y-1)^p={t'}^m$, o\`u $t'$ est une uniformisante de $\bar A$. Il en 
r\'esulte que $\bar y-1$ est la puissance $m$-i\`eme d'une uniformisante de 
$\bar B$. Comme $B$ est 
hens\'elien, $y=1+Z^m$, o\`u $Z$ est un param\`etre de $B$. Le calcul de 
$\sigma Z$ est imm\'ediat.
\par
(B) On suppose \`a pr\'esent que $0<\delta<v_K(p)$, alors, d'apr\`es la 
proposition \ref{reducmuptorseur}, $\spec B\to \spec A$ est un 
$fppf$-torseur sous $\mathcal H_n$, o\`u $v_K(p)=\delta+(p-1)n$. Il est donc 
donn\'e par
$$B=\frac{A[w]}{(\frac{(\pi^nw+1)^p-1}{\pi^{pn}}-u)}.$$
Il r\'esulte de la proposition \ref{reducmuptorseur} que $\bar u$ est unique 
\`a l'addition d'une puissance $p$-i\`eme pr\`es, donc la diff\'erentielle 
$d\bar u$ est bien d\'etermin\'ee. Les assertions suivantes se d\'emontrent 
de fa\c con analogue au cas multiplicatif.
\end{proof}
\begin{defn}
\label{diffassociee}
\par
Soit $\phi: \spec B \to \spec A$ un torseur comme dans le corollaire 
pr\'ec\'edent, on notera $\omega_{\phi}$ la $1$-forme diff\'erentielle 
\'egale \`a $\displaystyle{\frac{d\bar u}{\bar u}}$ si 
la r\'eduction est multiplicative, donn\'ee par l'\'equation $y^p=u$, et 
\'egale \`a $d\bar u$ si la r\'eduction est additive ou \'etale, donn\'ee par
l'\'equation $(\pi^nw+1)^p=1+\pi^{pn}u$. On dira que $\omega_{\phi}$ est 
la $1$-forme diff\'erentielle associ\'ee au torseur $\phi$.
\end{defn}
\begin{prop}\ \\
\label{pentediff}
\end{prop}
\begin{proof}
\par
Quitte \`a \'echanger les bords, on peut supposer $m\ge 0$. Comme 
$\varphi_K$ est \'etale, 
la fonction $\displaystyle{\frac{\sigma Z}{Z}}-1$ ne poss\`ede pas de z\'ero 
sur la couronne, donc d'apr\`es le lemme 1.6 de \cite{He2}, on a 
alors $m+m'=0$. De plus, on peut \'ecrire 
$\displaystyle{\frac{\sigma(Z)}{Z}-1}=\pi^{q}Z^mU,$
o\`u $q$ est un entier et $U$ est un inversible de 
$\mathcal O(\mathcal C_e)$. Alors, 
$$d_{\rho}=(p-1)v_{\rho}\left(\displaystyle{\frac{\sigma(Z)}{Z}-1}\right)=
q(p-1)+m(p-1)v_K(\rho).$$
 On conclut en remarquant que 
$q(p-1)=(p-1)v_{\eta}\left(\displaystyle{\frac{\sigma(Z)}{Z}-1}\right)
=d_\eta$. 
Enfin, la derni\`ere assertion r\'esulte du corollaire \ref{reducbord}.
\end{proof}
\section{Arbres de Hurwitz}
\par
Soit $\Gamma$ un arbre fini connexe, orient\'e \`a partir d'un sommet fix\'e 
$r_o$ appel\'e racine de l'arbre. 
On notera $\som \Gamma$ (resp. $\aret \Gamma$) l'ensemble des sommets 
(resp. ar\^etes) de $\Gamma$. Pour tout sommet $s$ de $\Gamma$, $\aret(s)$ 
(resp. $\aret^+(s)$)
sera l'ensemble des ar\^etes (resp. ar\^etes positives) d'origine $s$ ; le 
nombre d'ar\^etes d'origine $s$ est la valence du sommet $s$. 
Pour toute ar\^ete $a$ de $\Gamma$, 
$o_\Gamma(a)$, $t_\Gamma(a)$ et $\bar a$ d\'esignent respectivement 
l'origine, le sommet terminal et l'ar\^ete oppos\'ee de $a$. 
Une $n$-cha\^\i ne dans $\Gamma$ 
est un sous-arbre $\Gamma_0$ avec $\som \Gamma_0=\{s_i, 0 \le i \le n\}$ et 
$\aret \Gamma_0=\{a_i, \bar a_i, 0 \le i < n\}$, o\`u $o_{\Gamma_0}(a_i)=s_i$ 
et $t_{\Gamma_0}(a_i)=s_{i+1}$. La cha\^\i ne $\Gamma_0$ est 
positive si $a_i$ est positive pour tout $i$. Pour tout sommet $s$ il existe 
une unique chaine d'origine $r_o$ et de sommet terminal $s$, on la notera 
$\Gamma_s$.
\par
L'orientation de l'arbre induit une relation 
d'ordre sur $\som \Gamma$, \`a savoir $s_1 \le s_2$ si et seulement si 
il existe une 
cha\^ \i ne positive dans $\Gamma$ d'origine $s_1$ et d'extr\'emit\'e $s_2$. 
Un sommet maximal sera un sommet maximal pour cette relation d'ordre.
\begin{defn}
Etant donn\'e un arbre $\Gamma$, fini et connexe, on consid\`ere les 
donn\'ees suivantes :
\begin{itemize}
\item
Un sommet $r_0$ de $\Gamma$ ; on oriente $\Gamma$ \`a partir de $r_0$.
\item
Un entier $d_0$ divisible par $p-1$, v\'erifiant $0 \le d_0 \le v_K(p)$.
\item
Une application $\epsilon$ de $\aret \Gamma$ dans 
$\mathbb N$ telle que pour toute ar\^ete $a$ on ait 
$\epsilon(a)=\epsilon(\bar a)$. En particulier, $\epsilon$ d\'efinit une 
m\'etrique sur le graphe obtenu en retirant les ar\^etes $a$ telles que 
$\epsilon(a)=0$.
\item
Une application $m$ de $\aret \Gamma$ dans $\mathbb Z$.
\item
Une application $h$ de $\aret \Gamma$ dans $\mathbb Z/p\mathbb Z$.
\end{itemize}
\par
On d\'efinit alors, pour tout sommet $s$, un entier $d(s)$ par :
$$d(s):=d_0+(p-1)(\sum_{a \in \aret^+ \Gamma_s}m(a)\epsilon(a)).$$ 
On dira que la donn\'ee $\mathcal H:=(r_0, d_0, \epsilon, m, h)$ est 
une {\bf donn\'ee de Hurwitz} sur $\Gamma$, ou encore que 
$(\Gamma, \mathcal H)$ est un arbre de Hurwitz 
(sous-entendu d\'efini sur $K$), si les conditions 
$H[i]$, $1 \le i \le 7$, ci-dessous sont r\'ealis\'ees :
\newline
$H[1]$ Pour toute ar\^ete $a$, $m(\bar a)=-m(a)$ et $h(\bar a)=-h(a)$.
\newline
$H[2]$ Pour toute ar\^ete $a$, 
$m(a)=0$ si et seulement si $h(a) \ne 0 \mod p$. De plus, si $m(a) \ne 0$, 
$m(a)$ est premier \`a $p$.
\newline
$H[3]$ Un sommet $s$ de valence strictement sup\'erieur \`a $1$ 
est l'origine d'au moins trois ar\^etes distinctes, et on a les relations :
$$\sum_{a \in ar(s)}(m(a)+1)=2 \qquad et \qquad 
\sum_{a \in ar(s)}h(a)=0 \mod p.$$
$H[4]$ Si $a$ est une ar\^ete positive telle que $\epsilon(a)=0$, alors son 
sommet terminal est maximal. Si de plus $m(a)=0$, on dira 
que $a$ est une {\bf feuille}.
\newline
$H[5]$ Pour tout sommet $s$ de $\Gamma$, $0 \le d(s) \le v_K(p)$.\\
\par
On dira alors qu'un sommet $s$ de $\Gamma$ est :
\begin{itemize}
\item
{\bf  multiplicatif} si $d(s)=v_K(p)$,
\item
{\bf additif} si $0<d(s)<v_K(p)$,
\item
{\bf \'etale} si $d(s)=0$,
\end{itemize}
$H[6]$ Si $a$ est une feuille, l'origine de $a$ est un 
sommet multiplicatif.
\newline
$H[7]$ Si $s$ est additif ou \'etale, pour toute ar\^ete $a$ d'origine 
$s$, $h(a)=0 \mod p$.
\end{defn}
\par 
Si $(\Gamma_1, \mathcal H_1)$ et $(\Gamma_2, \mathcal H_2)$ sont 
des arbres de Hurwitz, on dira qu'ils sont \'equivalents si il existe un 
isomorphisme de $\Gamma_1$ sur $\Gamma_2$ qui transporte $\mathcal H_1$ sur 
$\mathcal H_2$.

\begin{defn}
Soit $(\Gamma,\mathcal H)$ un arbre de Hurwitz, et $a$ une ar\^ete positive 
de $\Gamma$. Consid\'erons le sous-arbre $\Gamma[a]$ dont 
l'ensemble des sommets est constitu\'e de l'origine de $a$ et des 
sommets $s$ de $\Gamma$ tels que $t_\Gamma(a)\le s$. La donn\'ee de Hurwitz 
$\mathcal H$ induit alors une donn\'ee de Hurwitz 
$\mathcal H[a]:=(o_\Gamma(a),d(o_\Gamma(a)),\epsilon,m,h)$ 
sur $\Gamma[a]$.
\end{defn}
\begin{defn}\label{dfnrealisable}
Soit $(\Gamma,\mathcal H)$ un arbre de Hurwitz.
\newline\newline
$[M]$ Si $s$ est un sommet de 
$\Gamma$ multiplicatif de valence sup\'erieure ou \'egale \`a $3$, on dira 
qu'il est r\'ealisable si il existe une fraction rationnelle 
$\bar u_s$ dans $k(t)$ et une injection $j_s$ de $\aret (s)$ dans 
$\mathbb P^1(k)$ telles que :
\begin{itemize}
\item
Si $a$ est une ar\^ete d'origine $s$, le r\'esidu de $\displaystyle
{\frac{d\bar u_s}{\bar u_s}}$ en $j_s(a)$ est \'egal \`a $h(a)$.
\item
$\Div(\displaystyle
{\frac{d\bar u_s}{\bar u_s}})=-\sum_{a \in \aret(s)}(m(a)+1)[j_s(a)].$
\end{itemize}
$[A]$ Si $s$ est un sommet de $\Gamma$ additif  de valence sup\'erieure ou 
\'egale \`a $3$, on dira qu'il est 
r\'ealisable si il existe une fraction rationnelle 
$\bar u_s$ dans $k(t)$ et une injection $j_s$ de $\aret (s)$ dans 
$\mathbb P^1(k)$ telles que :
$$\Div(d\bar u_s)=-\sum_{a \in \aret(s)}(m(a)+1)[j_s(a)].$$
\end{defn}
\begin{rem}
Tr\`es souvent, nous aurons \`a consid\'erer le cas o\`u un seul des 
$m(a)$ est positif (par exemple pour l'\'etude du disque formel). Le 
probl\`eme de l'existence d'une forme diff\'erentielle comme ci-dessus 
se ram\`ene \`a l'existence d'un point \`a coordonn\'ees distinctes d'une 
certaine sous-vari\'et\'e alg\'ebrique, en g\'en\'eral non irr\'eductible, 
d'un espace affine sur $k$ (voir la preuve de \ref{critrealdisc}). La 
dimension de l'espace de ces formes diff\'erentielles solutions de notre 
probl\`eme est alors $\displaystyle{[\frac mp]}$, o\`u $m$ est l'unique 
$m(a)$ positif (voir la remarque \ref{dimension}). En particulier, nous 
retrouvons la finitude de l'ensemble des solutions dans le cas o\`u $m<p$ 
(\cite{G-M 2}).
\end{rem}
\section{Automorphismes d'ordre $p$ du disque formel}
\subsection{Arbre de Hurwitz associ\'e \`a un automorphisme 
d'ordre $p$ du disque formel}
\par
Soit $\sigma$ un $R$-automorphisme d'ordre $p$ de $R[[Z]]$ (resp. $R\{Z\}$), 
on note $F_\sigma$ l'ensemble des 
points fixes g\'eom\'etriques de $\sigma$ dans la fibre g\'en\'erique 
$\mathcal D_K$ (resp. $D_K$) de $\mathcal D:=\spec R[[Z]]$ 
(resp. $D:=\spec R\{Z\}$); on 
supposera que le cardinal de $F_\sigma$ est \'egal \`a $m+1$ et que 
les points de $F_\sigma$ sont tous rationnels 
sur $K$ (le cardinal de $F_\sigma$ \'etant fini, on peut toujours le 
supposer, quitte \`a faire une extension finie de $K$). 
On consid\`ere le mod\`ele minimal 
$\mathcal M(\mathcal D_K,F_{\sigma})$ (resp. $\mathcal M(D_K,F_{\sigma})$) 
d\'eployant les sp\'ecialisations 
des points fixes en des points lisses et distincts, 
soit $\Gamma_\sigma^*$ l'arbre dual de sa fibre 
sp\'eciale, orient\'e \`a partir du sommet $r_0$ correspondant \`a la 
transform\'ee stricte du point g\'en\'erique 
de $\mathcal D_k$ (resp. de la fibre sp\'eciale de $D$).
Pour tout sommet $s$ de $\Gamma_\sigma^*$, notons $F_s$ le sous-ensemble, 
\'eventuellement vide, de 
$F_\sigma$ form\'e des points qui se sp\'ecialisent sur la composante 
correspondant \`a $s$. On note $\Gamma_\sigma$ l'arbre orient\'e 
d\'efini par :
\begin{itemize} 
\item
L'ensemble des sommets de $\Gamma_\sigma$ est la r\'eunion disjointe 
de $\som \Gamma_\sigma^*$ et de $F_\sigma$ (resp. de $\som \Gamma_\sigma^*$, 
de $F_\sigma$ et d'un singleton $\{s_\infty\}$). Si $x$ est un point 
de $F_\sigma$, on notera $s_x$ le sommet associ\'e.
\item
L'ensemble des ar\^etes positives de $\Gamma_\sigma$ est la r\'eunion 
disjointe de $\aret^+ \Gamma_\sigma^*$ et de $F_\sigma$ 
(resp. de $\som \Gamma_\sigma^*$, 
de $F_\sigma$ et d'un singleton $\{a_\infty\}$). Si $x$ est un point 
de $F_\sigma$, on notera $a_x$ l'ar\^ete positive associ\'ee.
\item
Si $a$ est une ar\^ete positive de $\Gamma_\sigma^*$, on a 
$o_{\Gamma_\sigma}(a)=o_{\Gamma_\sigma^*}(a)$ et 
$t_{\Gamma_\sigma}(a)=t_{\Gamma_\sigma^*}(a)$ ; 
pour tout point $x$ de $F_\sigma$, si $s$ est l'unique sommet de 
$\Gamma_\sigma^*$ tel que $x \in F_s$, alors 
$o_{\Gamma_\sigma}(a_x)=s$ et $t_{\Gamma_\sigma}(a_x)=s_x$ ; enfin, dans le 
cas du disque ferm\'e formel, on a $o_{\Gamma_\sigma}(a_\infty)=r_0$ et 
$t_{\Gamma_\sigma}(a_\infty)=s_\infty$. 
\end{itemize}
On d\'efinit une donn\'ee de Hurwitz 
$\mathcal H_\sigma:=(r_0,d_0,\epsilon,m,h)$ sur $\Gamma_{\sigma}$ 
de la fa\c con suivante :
\begin{itemize}
\item
L'entier $d_0$ est \'egal \`a la valuation de la diff\'erente au bord $\eta$ du 
disque formel.
\item
Si $a$ est une ar\^ete de $\Gamma_{\sigma}^*$, $\epsilon (a)$ est \'egal \`a 
l'\'epaisseur du point double correspondant dans la fibre sp\'eciale de 
$\mathcal M(\mathcal D_K,F_{\sigma})$ ; sinon, on pose $\epsilon(a)=0$.
\item
Soit $a$ une ar\^ete de $\Gamma_{\sigma}$. 
Si $a$ est une ar\^ete de $\Gamma_{\sigma}^*$ (resp. $a=a_x$ o\`u $x$ est un 
point de $F_{\sigma}$), elle correspond \`a un 
point double orient\'e $z_a$ de 
$\mathcal M_\sigma:=\mathcal M(\mathcal D_K,\sigma)$ (resp. \`a la 
sp\'ecialisation $\bar x$ de $x$ dans la fibre sp\'eciale de 
$\mathcal M_{\sigma}$). 
Soit $\xi_a$ le bord de la couronne formelle 
$\spec \hat {\mathcal O_{\mathcal M_\sigma,z_a}}$ (resp. le bord du disque formel 
$\spec \hat {\mathcal O_{\mathcal M_\sigma,\bar x}}$)
$\spec \hat {\mathcal O_{\mathcal M_\sigma,z_a}}$ qui correspond 
\`a l'origine de $a$. On notera $\hat {\mathcal O_a}$ le 
localis\'e-compl\'et\'e 
$(\hat {\mathcal O_{\mathcal M_\sigma,z_a})_{\xi_a}}^\wedge$, et $\omega_a$ 
la $1$-forme diff\'erentielle associ\'ee au torseur 
$\spec \hat {\mathcal O_a} \to \spec \hat {\mathcal O_a^{\sigma}}$ s'il a 
r\'eduction radicielle (voir \ref{diffassociee}). On pose alors 
$m(a):=-(1+\ord_{\xi_a}\omega_a)$ si ce torseur a r\'eduction radicielle, 
et sinon $m(a)$ est le conducteur de Hasse de l'extension 
$\hat {\mathcal O_a}\otimes_R k/\hat {\mathcal O_a^{\sigma}}\otimes_R k$.
On pose de plus $h(a)=0$ si la r\'eduction du torseur est additive ou 
\'etale et $h(a)=\textrm{R\'es}_{\xi_a}\omega_a$ si la r\'eduction est 
multiplicative. Par ailleurs, $m(\bar a_x)=-m(a_x)$ et $h(\bar a_x)=-h(a_x)$.
 Enfin, dans le cas du disque ferm\'e formel, on note 
$\hat {\mathcal O_\infty}:=R[[Z^{-1}]]\{Z\}$, et $\omega_\infty$ la $1$-forme 
diff\'erentielle associ\'ee au torseur 
$\spec \hat {\mathcal O_\infty} \to \spec \hat {\mathcal O_\infty^{\sigma}}$. 
On pose alors 
\text{$m(a_\infty)=-(1+\ord_{\infty}\omega_\infty)$}, $h(a_\infty)=0$ si la 
r\'eduction du torseur est additive ou \'etale et 
$m(a_\infty)=-(1+\ord_{\infty}\omega_\infty)$, 
$h(a_\infty)=\textrm{R\'es}_{\infty}\omega_\infty$ si la r\'eduction est 
multiplicative.
\end{itemize}
\begin{prop}
Avec les d\'efinitions ci-dessus, $(\Gamma_\sigma,\mathcal H_\sigma)$ est un 
arbre de Hurwitz. De plus :
\begin{itemize}
\item
Pour toute ar\^ete positive $a$, $m(a) \ge 0$.
\item
Si $a$ est une feuille, son sommet origine est un sommet maximal de 
$\Gamma_{\sigma}^*$. En d'autres termes, les sp\'ecialisations des 
points de $F_{\sigma}$ sont dans les bouts de la fibre sp\'eciale de 
$\mathcal M_{\sigma}$. Il en r\'esulte que pour toute ar\^ete positive $a$, en 
notant $F_a$ l'ensemble des feuilles de $\Gamma[a]$, 
$m(a)+1=\card(F_a)$.
\item
Tout sommet de valence sup\'erieure \`a $3$ est r\'ealisable.
\end{itemize}
\end{prop}
\begin{rem}
Cette proposition est une reformulation, dans le langage des arbres de 
Hurwitz, des r\'esultats de \cite{G-M 2}. La preuve ci-dessous est \`a 
quelques d\'etails pr\`es la m\^eme que celle expos\'ee dans l'article 
cit\'e au-dessus, nous ne la redonnons que pour illustrer la notion d'arbre de 
Hurwitz. On remarquera que dans le cas des arbres de Hurwitz provenant d'un 
automorphisme d'ordre $p$ d'un disque formel, la donn\'ee de l'arbre et de la 
racine permet de retrouver la valeur de $m(a)$ pour toute ar\^ete $a$. Ce ne 
sera plus le cas pour les arbres de Hurwitz provenant d'un automorphisme 
d'une couronne formelle.
\end{rem}
\begin{defn}
On dira qu'un arbre de Hurwitz provient d'un automorphisme $\sigma$ 
d'ordre $p$ du disque formel (ouvert ou ferm\'e) s'il est \'equivalent \`a 
$(\Gamma_{\sigma}, \mathcal H_{\sigma})$.
\end{defn}
\begin{proof}\ \\
$H[1]$ Si $a$ est une ar\^ete de $\Gamma_\sigma^*$, la relation 
$m(\bar a)=-m(a)$ r\'esulte de \ref{pentediff}. Sinon, elle d\'ecoule de 
la d\'efinition. Par ailleurs, la relation $h(\bar a)=-h(a)$ d\'ecoule 
\'egalement de la d\'efinition.
\newline\newline
$H[2]$ Soit $a$ une ar\^ete de $\Gamma_\sigma$, on peut, quitte \`a changer 
$a$ en $\bar a$, supposer $a$ positive. En particulier, l'origine de $a$ n'est 
pas un sommet maximal. Supposons tout d'abord l'origine de $a$ multiplicative, 
la forme diff\'erentielle $\omega_a$ 
est alors logarithmique, donc l'ordre de $\omega_a$ 
est sup\'erieur \`a $-1$. Alors $m(a)=0$ est 
\'equivalent \`a $\ord_{\xi_a}\omega_a=-1$, ce qui est \'equivalent \`a 
$h(a)=\textrm{R\'es}_{\xi_a}\omega_a \ne 0 \mod p$. Si maintenant on suppose 
l'origine de $a$ additive ou \'etale, on a par hypoth\`ese $h(a)=0$ et 
comme $\omega_a$ est maintenant une diff\'erentielle exacte, on ne peut pas 
avoir $\ord_{\xi_a}\omega_a=-1$. Ainsi, $m(a)\ne 0$.
\newline\newline
$H[3]$ Le fait qu'un sommet ne peut pas \^etre de valence \'egale \`a 
$2$ r\'esulte de la 
minimalit\'e du mod\`ele $\mathcal M_{\sigma}$. 
Soit $s$ un sommet de $\Gamma_{\sigma}^*$, correspondant \`a une droite 
projective $E_s$, et $\tilde E_s$ 
l'ouvert affine de $E_s$ compl\'ementaire des points doubles et des 
sp\'ecialisations des points de $F_s$. On note $E'_s$ la composante au-dessous 
de $E_s$ dans le quotient de $\mathcal M_{\sigma}$ par $\sigma$ et 
$\tilde E'_s$ l'image de $\tilde E_s$ dans $E'_s$.
\par
Soit $\hat {\mathcal M_{\sigma}}$ le compl\'et\'e formel de 
$\mathcal M_{\sigma}$ le long de sa fibre sp\'eciale, et 
$\tilde \mathcal E_s$ le sous-sch\'ema formel affine de 
$\hat {\mathcal M_{\sigma}}$ de fibre sp\'eciale $\tilde E_s$. On d\'efinit 
de mani\`ere analogue $\tilde \mathcal E'_s$. La $R$-alg\`ebre 
$A_s:=\mathcal O(\tilde \mathcal E'_s)$ est compl\`ete pour la topologie 
$\pi$-adique et factorielle. En notant $B_s:=\mathcal O(\tilde \mathcal E_s)$, 
on consid\`ere alors la r\'eduction du torseur 
$\spec B_s \to \spec A_s$, qui est radicielle.
\par
Si la r\'eduction est additive, donn\'ee par une \'equation 
$(\pi^nw+1)^p=1+\pi^{pn}u_s$, la diff\'erentielle $\omega_s:=d\bar u_s$ 
est uniquement d\'etermin\'ee. Comme $\tilde E_s$ est une courbe lisse, 
$\omega_s$ n'a ni z\'ero ni p\^ole sur $\tilde E'_s$. On a alors 

$$\Div {\omega}_s=-\sum_{a \in \aret(s)}(m(a)+1)[t_a].$$
(En notant $t_a$ l'image dans $E'_s$ du point de $E_s$ correspondant \`a 
l'ar\^ete $a$). La relation $\sum_{a \in \aret(s)}(m(a)+1)=2$ est alors 
simplement le fait que le degr\'e de $\omega_s$ est \'egal \`a $-2$. De plus, 
comme $s$ est additif, pour toute ar\^ete $a$ d'origine $s$, on a $h(a)=0$ 
par d\'efinition. Ainsi, la somme des $h(a)$ est clairement nulle. 
\par
Si la r\'eduction est multiplicative, donn\'ee par une \'equation 
$y^p=u_s$, la diff\'erentielle 
$\omega_s:=\displaystyle{\frac{d\bar u_s}{\bar u_s}}$ est uniquement 
d\'etermin\'ee. Comme $\tilde E_s$ est une courbe lisse, 
$\omega_s$ n'a ni z\'ero ni p\^ole sur $\tilde E'_s$. On a alors 
$$\Div \omega_s=-\sum_{a \in \aret(s)}(m(a)+1)[t_a].$$
et pour $a$ dans $\aret(s)$, le r\'esidu de $\omega_s$ en $t_a$ est $h(a)$. 
Les relations $\sum_{a \in \aret(s)}(m(a)+1)=2$ et 
$\sum_{a \in \aret(s)}h(a)=0$ traduisent alors respectivement que le degr\'e 
de $\omega_s$ est $-2$ et le th\'eor\`eme des r\'esidus.
\par
En particulier, on voit d\`es \`a pr\'esent que les sommets correspondant \`a 
des droites projectives sont r\'ealisables.
\newline\newline
$H[4]$ est imm\'ediat. (Les feuilles correspondent aux points de $F_{\sigma}$).
\newline\newline
$H[5]$ r\'esulte du 
\begin{lem}
\par
Soit $s$ un sommet de $\Gamma_\sigma^*$. Le point g\'en\'erique de la 
composante $E_s$ d\'efinit une valuation discr\`ete $v_s$ du corps des 
fractions $\mathcal K$ de $R[[Z]]$ (resp. $R\{Z\}$). 
L'entier $d(s):=d_0+(p-1)(\sum_{a \in \aret^+ \Gamma_s}m(a)\epsilon(a))$ 
est alors \'egal \`a la valuation de la diff\'erente de l'extension 
$(\mathcal K,v_s)/(\mathcal K^{\sigma},v_s)$.
\end{lem}
La preuve est imm\'ediate, en utilisant la proposition \ref{pentediff}, par 
r\'ecurrence sur le nombre de sommets de $\Gamma_s$.
\newline\newline
$H[6]$ Une feuille correspond \`a un point fixe de $\sigma$, en particulier, 
la composante qui le porte a pour diff\'erente $v_K(p)$, c'est-\`a-dire 
l'origine de la feuille est un sommet multiplicatif.
\newline\newline
$H[7]$ est imm\'ediat d'apr\`es la d\'efinition.
\par
La donn\'ee 
$\mathcal H_\sigma$ est donc bien une donn\'ee de Hurwitz sur $\Gamma_\sigma$. 
Soit $a$ une ar\^ete positive, et $s$ son origine. Consid\'erons le mod\`ele 
$M_a$ de $\mathcal D_K$ obtenu \`a partir de $\mathcal M_\sigma$ en 
recontractant les droites projectives correspondant aux sommets $s'$ 
tels que $\Gamma_{s'}$ passe par $a$. Ces droites se contractent en un point 
$z_a$ de $E_s$ qui est alors 
lisse sur $M_a$, soit $Z_a$ un param\`etre du disque formel 
$\hat {\mathcal O_{M_a,z_a}}$, i.e. $\hat {\mathcal O_{M_a,z_a}}=R[[Z_a]]$. 
On a alors, d'apr\`es la proposition \ref{pentediff}, 
$m(a)=\ord_{Z_a}\left( \displaystyle{\frac{\sigma Z_a}{Z_a}}-1\right) \ge 0$. 
Soit $a_0$ une feuille de $\Gamma_\sigma$. Son origine $s$ est un sommet 
de $\Gamma_\sigma^*$, multiplicatif. Supposons qu'il existe une ar\^ete $a$ 
positive de $\Gamma_\sigma^*$ de sommet origine $s$, on a alors 
$m(a)\ge 0$. Comme 
$d(t_{\Gamma_\sigma}(a))\le v_K(p)=d(o_{\Gamma_\sigma}(a))$, on doit 
avoir $m(a)=0$. Plus g\'en\'eralement, pour toute cha\^ \i ne positive 
d'origine $s$, le sommet terminal est multiplicatif. Prenons une telle 
cha\^ \i ne, maximale. La derni\`ere ar\^ete est une feuille. De son 
sommet origine $s'$ ne partent que des feuilles par maximalit\'e. 
De la relation $\sum_{a' \in \aret(s')}(m(a')+1)=2$ et du fait qu'il 
existe au moins deux feuilles d'origine $s'$, on d\'eduit la contradiction 
$m(a)<0$ pour l'unique ar\^ete positive de sommet terminal $s'$. Ainsi, toute 
ar\^ete d'origine $s$ est une feuille.
\par
Soit $a$ une ar\^ete positive de $\Gamma_\sigma$. Si $a$ est une feuille, 
alors $m(a)=0$, et \text{$\card F_a=1$.} Si $a$ n'est pas une feuille, le 
sommet 
terminal $s$ de $a$ est de valence sup\'erieure \`a $3$. Donc, d'apr\`es 
$H[3]$, $m(a)+1=\sum_{a' \in \aret^+(s)}(m(a')+1)$. Par r\'ecurrence sur le 
nombre maximal de sommets d'une cha\^ \i ne positive reliant le sommet 
terminal de $a$ \`a une feuille, on peut alors supposer 
$m(a')+1=\card F_{a'}$ pour toute ar\^ete positive $a'$ d'origine $s$. Comme 
$F_a$ est la r\'eunion disjointe des $F_{a'}$, on a alors 
\text{$m(a)+1=\card F_a$.} 
Remarquons qu'on peut aussi d\'eduire cette formule du th\'eor\`eme de 
pr\'eparation de Weierstrass.
\end{proof}
\subsection{Th\'eor\`eme de r\'ealisation pour le disque formel}
\begin{thm}\ \\
\label{hurwitzdisc}
\par
Soit $(\Gamma, \mathcal H)$ un arbre de Hurwitz, avec 
$\mathcal H=(r_0,d_0,\epsilon,m,h)$ ; en notant $R$ l'anneau de valuation 
de $K$, $(\Gamma, \mathcal H)$ provient d'un $R$-automorphisme d'ordre $p$ 
de $R[[Z]]$ si et seulement si il v\'erifie les conditions $D[i]$ ci-dessous 
($1 \le i \le 3$) :
\newline
$D[1]$ La racine $r_0$ est de valence $1$.
\newline
$D[2]$ Toute ar\^ete aboutissant \`a un sommet maximal est une feuille.
\newline
$D[3]$ Tout sommet de valence sup\'erieure ou \'egale \`a $3$ est r\'ealisable.
\end{thm}
\par
La partie directe du th\'eor\`eme, due \`a B. Green et M. Matignon, 
a d\'ej\`a \'et\'e prouv\'ee, \`a 
l'exception de $D[1]$, qui r\'esulte de la construction du mod\`ele minimal : 
En fait, le sommet terminal de l'unique ar\^ete d'origine $r_0$ correspond au 
plus petit disque ferm\'e contenant tous les points de $F_{\sigma}$.
\begin{rem}
Quitte \`a faire une extension finie de $K$, on peut toujours se ramener 
au cas o\`u la valuation $d_0$ de 
la diff\'erente au bord du disque formel est nulle, au moyen d'un recollement 
convenable d'une couronne formelle semi-ouverte avec action d'un 
automorphisme d'ordre $p$. En effet, soit $\sigma$ un automorphisme 
d'ordre $p$ du disque formel, avec $m+1$ points 
fixes dans la fibre g\'en\'erique. Posons $v_K(p)=d_0+(p-1)n_0$, d'apr\`es 
le corollaire \ref{reducbord}, il existe un param\`etre $Z$ 
au bord du disque formel tel que 
$$\sigma Z=\zeta^{-\frac 1m}Z(1+\lambda\pi^{-n_0}
\zeta^{-1}Z^{m})^{-\frac 1m}.$$
Consid\'erons maintenant l'automorphisme $\sigma_0$ d'ordre $p$ de la couronne 
formelle semi-ouverte 
$\spec \displaystyle{\frac{R[[V]]\{W\}}{(VW-\rho)}}$, donn\'e par :
$$\sigma_0 V:=\zeta^{-\frac 1m}V(1+V^{m})^{-\frac 1m},$$
avec $\rho^m:=\lambda\pi^{-n_0}\zeta^{-1}$. On utilise alors un lemme analogue au lemme \ref{patchingopendisc} pour construire un prolongement de 
$\sigma$ \`a un disque formel plus grand, en identifiant $W=Z^{-1}$, et le 
choix de $\sigma_0$ montre que 
l'automorphisme obtenu agit sans inertie au bord.
\end{rem}
\subsection{Construction d'automorphismes d'ordre $p$}
\par
Le but de cette partie est de donner les r\'esultats constituant le coeur 
technique de la preuve de la partie r\'eciproque du th\'eor\`eme 
\ref{hurwitzdisc}. Notamment, on va voir que 
l'\'etape essentielle est la construction 
d'automorphismes d'ordre $p$ du disque ferm\'e \`a partir d'un arbre de 
Hurwitz. 
On fixe une extension finie $K$ de $K_0$, et on note $R$ son anneau de 
valuation ; $\pi$ d\'esigne une uniformisante de $R$
\subsubsection{Lemmes de recollement}
\begin{lem}
\label{patchingcloseddisc}
On se donne des \'el\'ements $a_1,\dots,a_r$ de $R$ distincts deux \`a deux 
modulo $\pi$. 
On note, pour $i$ de $1$ \`a $r$, $\alpha_i$ (resp. $\beta_i$) l'injection 
canonique de la $R$-alg\`ebre 
$R\{Z,(Z-a_j)^{-1}\}_{1 \le j \le r}$ (resp. $R[[X_i]]$) dans 
$R[[Z-a_i]]\{(Z-a_i)^{-1}\}$ (resp. $R[[X_i]]\{X_i^{-1}\}$). On se donne 
\'egalement des isomorphismes de $R$-alg\`ebres 
$\psi_i: R[[X_i]]\{X_i^{-1}\} \to R[[Z-a_i]]\{(Z-a_i)^{-1}\}$.
\par
Si $\theta$ est l'homomorphisme de $R$-module de 
$R\{Z,(Z-a_j)^{-1}\}_{1 \le j \le r}\times \prod_{1 \le i \le r}R[[X_i]]$ 
dans $\prod_{1 \le i \le r}R[[Z-a_i]]\{(Z-a_i)^{-1}\}$ d\'efini par 
$$\theta(f_0,f_1,\dots,f_r)=(\alpha_1(f_0)-\psi_1 \circ \beta_1(f_1),\dots,
\alpha_r(f_0)-\psi_r \circ \beta_r(f_r)),$$ 
pour $f_0$ dans $R\{Z,(Z-a_j)^{-1}\}_{1 \le j \le r}$ et $f_i$ dans 
$R[[X_i]]$ pour $1 \le i \le r$, alors $\theta$ est surjective 
et son noyau $N$ est une $R$-alg\`ebre de s\'eries enti\`eres restreintes en 
une variable (autrement dit de la forme $R\{Z_0\}$).
\end{lem}
\begin{proof}
\par
On voit imm\'ediatement que $N$ est une sous-$R$-alg\`ebre de 
$$R\{Z,(Z-a_j)^{-1}\}_{1 \le j \le r}\times \prod_{1 \le i \le r}R[[X_i]].$$
On posera $z:=Z \mod \pi$ et, pour 
$i$ de $1$ \`a $r$, $x_i:=X_i \mod \pi$. De l'homomorphisme de $R$-modules 
$\theta$, on obtient un homomorphisme de $k$-espace vectoriel 
$$\bar \theta:=\theta\otimes_Rk : k[z,(z-\bar a_j)^{-1}]_{1 \le j \le r} 
\times \prod_{1 \le i \le r} k[[x_i]] \to \prod_{1 \le i \le r} k((x_i)).$$
\par
Comme $\bar \psi_i(\bar \alpha_i(x_i))$ est une uniformisante de 
$k((z-\bar a_i))$, quitte \`a changer le param\`etre $X_i$, on peut 
supposer que $\bar \psi_i(\bar \alpha_i(x_i))=z-\bar a_i$. On voit alors 
imm\'ediatement que $\bar \theta$ est surjectif et son noyau est \'egal 
\`a $k[z]$. Il suit alors que $\theta$ est surjectif, et on obtient donc la 
suite exacte :
$$0 \to N \to 
R\{Z,(Z-a_j)^{-1}\}_{1 \le j \le r}\times \prod_{1 \le i \le r}R[[X_i]]
\to \prod_{1 \le i \le r}R[[Z-a_i]]\{(Z-a_i)^{-1}\} \to 0$$
Alors, $N\otimes_Rk$ 
s'identifie au noyau $k[z]$ de $\bar \theta$ car 
$\prod_{1 \le i \le r}R[[Z-a_i]]\{(Z-a_i)^{-1}\}$ est plat sur $R$. 
Relevons $z$ en un \'el\'ement $Z_0$ de $N$, on a alors $N=R\{Z_0\}$.
\end{proof}
On d\'emontre de mani\`ere identique le lemme :
\begin{lem}\label{patchingopendisc}
Soient $e$ une entier strictement positif, $\beta$ l'injection canonique 
de $\displaystyle{\frac{R[[Z_1,Z_2]]}{(Z_1Z_2-\pi^e)}}$ dans 
$R[[Z_2]]\{Z_2^{-1}\}$, $\alpha$ l'injection canonique de $R\{Z\}$ dans 
$R[[Z^{-1}]]\{Z\}$ et $\psi$ un isomorphisme de $R$-alg\`ebres de 
$R[[Z_2]]\{Z_2^{-1}\}$ sur $R[[Z^{-1}]]\{Z\}$. Si $\theta$ d\'esigne 
l'homomorphisme de $R$-modules 
$$\theta : R\{Z\} \times \frac{R[[Z_1,Z_2]]}{(Z_1Z_2-\pi^e)} \to 
R[[Z^{-1}]]\{Z\}$$
d\'efini par $\theta(f,g):=\alpha(f)-\psi \circ \beta (g)$, alors $\theta$ 
est surjectif, et son noyau $N$ est une alg\`ebre de s\'eries enti\`eres en 
une variable, \`a coefficients dans $R$ (autrement dit de la forme $R[[Z_0]]$).
\end{lem}
\subsubsection{Application aux arbres de Hurwitz}
\begin{prop}
\label{hfmult}
Soit $\Gamma$ un arbre fini connexe et 
$\mathcal H=(r_1,d_1,\epsilon,m,h)$ 
une donn\'ee de Hurwitz sur $\Gamma$ ; on suppose que $(\Gamma, \mathcal H)$ 
v\'erifie les conditions suivantes :
\begin{itemize}
\item
La racine $r_1$ est de valence sup\'erieure ou \'egale \`a $3$
\item
Toute ar\^ete $a$ d'origine $r_1$ v\'erifie 
$\epsilon(a)=0$. En particulier, le seul sommet non maximal est 
la racine $r_1$.
\item
Il existe une unique ar\^ete $a_1$ d'origine $r_1$ qui n'est pas une feuille 
(i.e. $m(a_1) \ne 0$).
\item
La racine $r_1$ est r\'ealisable.
\end{itemize}
Alors, si $R$ d\'esigne l'anneau de valuation de $K$, 
il existe un automorphisme $\sigma$ d'ordre $p$ du disque ferm\'e 
formel $\spec R\{Z\}$ dont l'arbre de Hurwitz associ\'e est \'equivalent \`a 
$(\Gamma, \mathcal H)$, et tel qu'il existe un param\`etre $V$ de 
$R[[Z^{-1}]]\{Z\}$ avec 
$$\sigma(V)=\zeta^{-\frac 1{m(a_1)}} 
V(1 +\pi^{v_K(\lambda)}V^{m(a_1)})^{-\frac 1{m(a_1)}}.$$
\end{prop}
\begin{figure}[htpb]
         \input{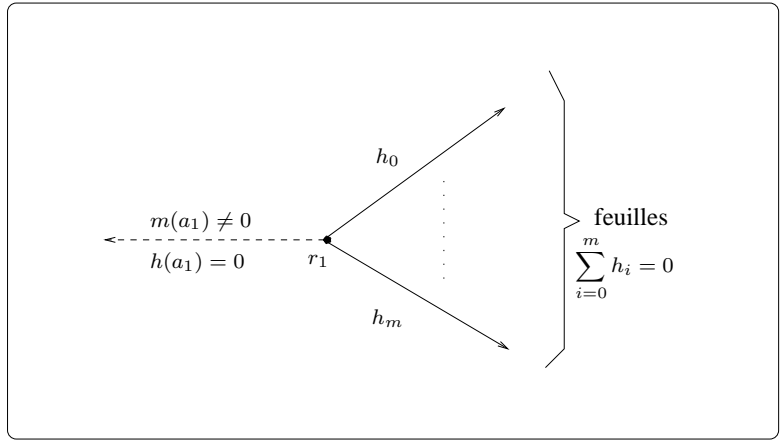}
         \caption{Premi\`ere r\'ealisation (\ref{hfmult})}
\end{figure}
\begin{proof}
\par
Remarquons tout d'abord que si le nombre de feuilles $a$ de 
$(\Gamma, \mathcal H)$ est \text{$m+1$,} alors on a $m(a_1)=-m$ 
d'apr\`es la condition $H[3]$. De plus, la racine est multiplicative 
d'apr\`es $H[6]$. 
Soit $\bar u_{r_1}$ et $j_{r_1}$ v\'erifiant les conditions de $[M]$. 
Pour simplifier les notations, on notera $\bar u:=\bar u_{r_1}$, et 
$t_a:=j_s(a)$ pour toute ar\^ete $a$ d'origine $r_1$. On peut supposer 
$t_{a_1}=+\infty$. Alors, quitte \`a multiplier $\bar u$ par une puissance 
$p$-i\`eme, ce qui ne change pas la diff\'erentielle 
$\displaystyle{\frac{d\bar u}{\bar u}}$, on peut supposer que 
$\bar u=\prod_{a \in F}(t-t_a)^{h(a)}$, o\`u $F$ d\'esigne l'ensemble des 
feuilles distinctes de $a_1$. Soit $T_{a}$ 
dans $R$ relevant $t_{a}$, et $u$ un \'el\'ement de la $R$-alg\`ebre 
\text{$\mathcal A:=R\{T,(T-T_{a})^{-1}\}_{a \in F}$} relevant $\bar u$. 
Notons 
$\mathcal B:=\displaystyle{\frac {\mathcal A[Y]}{(Y^p-u)}}$, c'est une 
$\mathcal A$-alg\`ebre libre de rang $p$. On consid\`ere l'automorphisme 
de $\mathcal A$-alg\`ebre (d'ordre $p$) $\sigma$ tel que $\sigma(Y)=\zeta Y$. 
Le $k$-morphisme $\spec \bar \mathcal B \rightarrow \spec \bar \mathcal A$ 
est un rev\^etement purement ins\'eparable de degr\'e $p$ de l'ouvert 
affine $\mathbb P_k^1 \setminus \{\infty,t_a\}_{a \in F}$, et $d\bar u$ ne 
poss\`ede ni z\'ero, ni p\^ole sur cet ouvert, donc 
$\spec \bar \mathcal B$ est lisse, et c'est un ouvert de la droite projective 
sur $k$ de la forme $\mathbb P_k^1 \setminus \{\infty,z_a\}_{a \in F}$, o\`u 
les $z_a$ sont des 
\'el\'ements de $k$ distincts deux \`a deux. Autrement dit, on a  
$\bar \mathcal B=k[z,(z-z_a)^{-1}]_{a \in F}$. Soit $Z$ 
(resp. $Z_{a}$) un rel\`evement de $z$ (resp. $z_{a}$) dans 
$\mathcal B$ (resp. $R$), comme $\mathcal B$ est plat sur $R$, on a 
$\mathcal B=R\{Z,(Z-Z_{a})^{-1}\}_{a \in F}$. Soit $\omega$ la forme diff\'erentielle $\displaystyle{\frac{d\bar u}{\bar u}}$.
\par
Pour $a \in F$, le r\'esidu de $\omega$ en $t_a$ est $h(a)\ne 0 \mod p$ ; 
il r\'esulte 
alors de la proposition \ref{reducbord} (A,a) que 
$\mathcal B \otimes_{\mathcal A}R[[T-T_{a}]]\{(T-T_{a})^{-1})\}=
R[[Z-Z_{a}]]\{(Z-Z_{a})^{-1}\}$, et l'action induite par 
celle de $\sigma$ est donn\'ee sur un param\`etre convenable 
$X_{(a)}$ par $\sigma(X_{(a)})=\zeta^{h(a)^{-1}}(X_{(a)})$. Remarquons que 
cette action se prolonge canoniquement en une action sur $R[[X_{(a)}]]$. 
Soit $N$ le noyau de l'homomorphisme (surjectif) de 
$R[\langle \sigma \rangle]$-module $\phi:\mathcal B 
\times \prod_{a \in F}R[[X_{(a)}]] 
\rightarrow \prod_{a \in F}R[[Z-Z_{a}]]\{(Z-Z_{a})^{-1}\}$ 
d\'efini par $\phi(b,(f_{a}))=(b\otimes 1-f_{a})_{a \in F}$. 
D'apr\`es le lemme \ref{patchingcloseddisc}, la 
$R[\langle \sigma \rangle]$-alg\`ebre $N$ est alors de la forme $R\{Z_0\}$. 
De plus, l'arbre de Hurwitz associ\'e \`a l'automorphisme $\sigma$ de 
$R\{Z_0\}$ est isomorphe \`a $(\Gamma, \mathcal H)$. La derni\`ere assertion 
r\'esulte de la proposition \ref{reducbord} (A,b), car 
$\ord_{\infty}d\bar u=-(m(a_1)+1)$.
\end{proof}
\begin{cor}
\label{homult}
Soit $\Gamma$ un arbre fini connexe et $\mathcal H=(r_1,d_0,\epsilon,m,h)$ 
une donn\'ee de Hurwitz sur $\Gamma$ ; on suppose que $(\Gamma, \mathcal H)$ 
v\'erifie les conditions suivantes :
\begin{itemize}
\item
La racine $r_0$ est l'origine d'une unique ar\^ete $a_0$, et 
$\epsilon(a_0) \ne 0$, $m(a_0) \ne 0$.
\item
Si $a$ est une ar\^ete positive d'origine $r_1:=t_\Gamma(a_0)$, $a$ est une 
feuille. En particulier, $r_1$ est multiplicatif.
\item
Le sommet $r_1$ est r\'ealisable.
\end{itemize}
Alors, il existe un automorphisme $\sigma$ d'ordre $p$ du disque ouvert 
formel $\spec R[[Z]]$ dont l'arbre de Hurwitz associ\'e est \'equivalent \`a 
$(\Gamma, \mathcal H)$, et tel qu'il existe un param\`etre $X$ de 
$R[[Z]]\{Z^{-1}\}$ avec 
$$\sigma(X)=\zeta^{-\frac 1{m(a_0)}}X
(1+\pi^{v_K(\lambda )-n_0}X^{m(a_0)})^{-\frac 1{m(a_0)}},$$ 
o\`u $n_0$ est un entier 
donn\'e par la relation $v_K(p)=d_0+(p-1)n_0$.
\end{cor}
\begin{figure}[htpb]
         \input{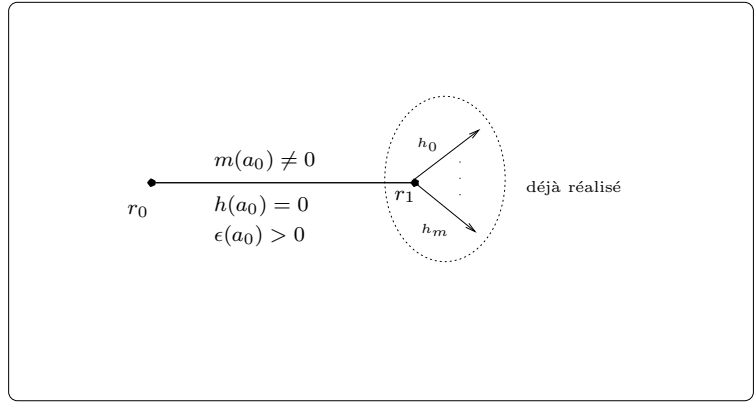}
         \caption{Deuxi\`eme r\'ealisation (\ref{homult})}
\end{figure}

\begin{proof}
\par
Soit $\mathcal H'$ la donn\'ee de Hurwitz sur $\Gamma$ d\'efinie par 
$$\mathcal H':=(r_1,d_1,\epsilon',m,h),$$ 
o\`u $d_1:=d_0+m(a_0)(p-1)\epsilon(a_0)$, 
$\epsilon'(a)=\epsilon(a)$ si $a \notin \{a_0,\bar a_0\}$, et 
$\epsilon'(a_0)=0$. Alors, $(\Gamma,\mathcal H')$ v\'erifie les hypoth\`eses 
de la proposition \ref{hfmult} (avec $a_1=\bar a_0$). Il existe ainsi un 
$R$-automorphisme $\sigma$ 
de $R\{Z\}$, d'ordre $p$, dont l'arbre de Hurwitz s'identifie \`a 
$(\Gamma,\mathcal H')$, et tel qu'il existe un param\`etre $V$ de 
$R[[Z^{-1}]]\{Z\}$ 
tel que l'action induite de $\sigma$ sur $V$ soit donn\'ee par 
$\sigma(V)=\zeta^{\frac 1m} V
(1+\pi^{v_K(\lambda)}V^{-m})^{\frac 1m}$, o\`u $m=m(a_0)>0$ et 
$\lambda=\zeta-1$. 
Notons $e:=\epsilon(a_0)$ et consid\'erons la $R$-alg\`ebre 
$\displaystyle{\frac {R[[Z_1,Z_2]]}{(Z_1Z_2-\pi^e)}}$, munie du 
$R$-automorphisme (d'ordre $p$) 
$\sigma(Z_1)=\zeta^{-\frac 1m} Z_1(1+\pi^{v_K(\lambda)-em}Z_1^m)^{-\frac 1m}$. 
Soit $\psi$ l'isomorphisme 
de $R$-alg\`ebres de $R[[Z_2]]\{Z_2^{-1}\}$ sur $R[[Z^{-1}]]\{Z\}$ qui envoie 
$Z_2$ sur $V$, $\psi$ est $\sigma$-\'equivariant ; en utilisant le lemme 
\ref{patchingopendisc}, 
on obtient un automorphisme $\sigma$ de $R[[Z_0]]$ dont l'arbre de Hurwitz 
associ\'e s'identifie \`a $(\Gamma, \mathcal H)$, et le param\`etre au bord 
$Z':=Z_1$ poss\`ede la propri\'et\'e attendue ; en effet, 
$v_K(p)=d(r_1)=d_0+em(a_0)(p-1)$, donc $n_0=em$.
\end{proof}
\begin{prop}
\label{hfadd}
Soit $\Gamma$ un arbre fini connexe et $\mathcal H=(r_1,d_1,\epsilon,m,h)$ 
une donn\'ee de Hurwitz sur $\Gamma$ ; on suppose que $(\Gamma, \mathcal H)$ 
v\'erifie les conditions suivantes :
\begin{itemize}
\item
La racine $r_1$ est additive.
\item
Il existe une unique ar\^ete $a_1$ d'origine $r_1$ 
v\'erifiant $m(a_1)<0$.
\item
Si $a$ est une ar\^ete d'origine $r_1$ distincte de $a_1$, l'arbre de 
Hurwitz $(\Gamma_a,\mathcal H_a)$ 
provient d'un $R$-automorphisme $\sigma_a$ d'ordre $p$ d'un disque formel 
$\mathcal D_a$, et il existe un param\`etre $X_a$ au bord de 
$\mathcal D_a$ tel que 
$$\sigma_{a}(X_a)=\zeta^{-\frac 1{m(a)}}X_a
(1+\pi^{v(\lambda)-n_1}X_a^{m(a)})^{-\frac 1{m(a)}},$$
o\`u 
$n_1$ est l'entier d\'efini par la relation $v_K(p)=d(r_1)+n_1(p-1)$.
\item
La racine $r_1$ est r\'ealisable.
\end{itemize}
Alors, si $R$ d\'esigne l'anneau de valuation de $K$, 
il existe un automorphisme $\sigma$ d'ordre $p$ du disque ferm\'e 
formel $\spec R\{Z\}$ dont l'arbre de Hurwitz associ\'e est isomorphe \`a 
$(\Gamma, \mathcal H)$, et tel qu'il existe un param\`etre $V$ de 
$R[[Z^{-1}]]\{Z\}$ avec 
$$\sigma(V)=\zeta^{-\frac 1{m(a_1)}} 
V(1 +\pi^{v_K(\lambda)-n_1}V^{m(a_1)})^{-\frac1{m(a_1)}}.$$
\end{prop}
\begin{figure}[htpb]
         \input{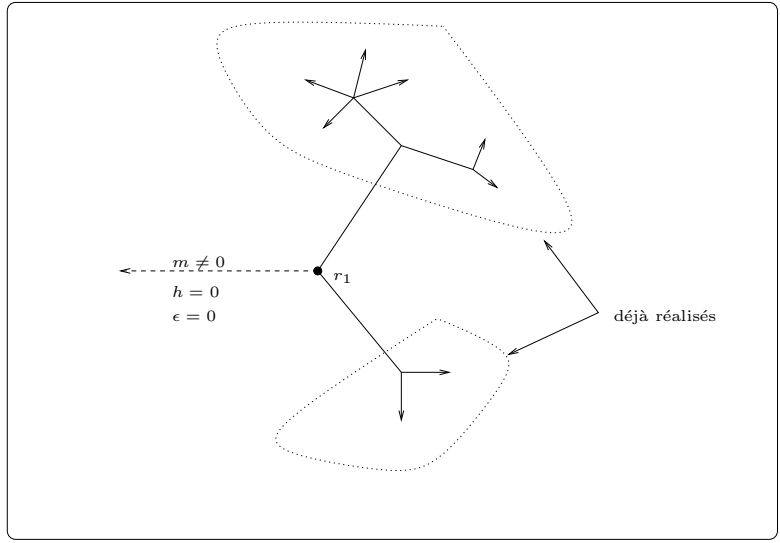}
         \caption{Troisi\`eme r\'ealisation (\ref{hfadd})}
\end{figure}
\begin{proof}
\par
Par la condition $H[3]$, on a 
$$m:=m(\bar a_1)=-m(a_1)=-1+\displaystyle{\sum_{a \in \aret (r_1)\setminus \{a_1\}}(m(a)+1)}.$$
Soit $\bar u_{r_1}$ et $j_{r_1}$ v\'erifiant les conditions de $[A]$. 
Pour simplifier les notations, on notera $\bar u:=\bar u_{r_1}$, et 
$t_a:=j_s(a)$ pour toute ar\^ete $a$ d'origine $r_1$. On peut supposer 
$t_{a_1}=+\infty$. Alors, quitte \`a ajouter \`a $\bar u$ par une puissance 
$p$-i\`eme, ce qui ne change pas la diff\'erentielle 
$d\bar u$, on peut supposer que les p\^oles de $\bar u$ sont les $t_a$, pour 
$a$ dans $\aret (r_1)\setminus \{a_1\}$. Soit $T_{a}$ 
dans $R$ relevant $t_{a}$, et $u$ un \'el\'ement de la $R$-alg\`ebre 
$\mathcal A:=R\{T,(T-T_{a})^{-1}\}_{a \in F}$ relevant $\bar u$. 
Notons 
$\mathcal B:=\displaystyle{\frac {\mathcal A[W]}{(\frac{(\pi^{n_1}W+1)^p-1}
{\pi^{pn_1}}-u)}}$, c'est une 
$\mathcal A$-alg\`ebre libre de rang $p$. On consid\`ere l'automorphisme 
de $\mathcal A$-alg\`ebre (d'ordre $p$) $\sigma$ tel que 
$\sigma(W)=\zeta W+\frac{\zeta-1}{\pi^{n_1}}$. 
Le $k$-morphisme $\spec \bar \mathcal B \rightarrow \spec \bar \mathcal A$ 
est un rev\^etement purement ins\'eparable de degr\'e $p$ de l'ouvert 
affine $\mathbb P_k^1 \setminus \{\infty,t_a\}_{a \in \aret (r_1)\setminus \{a_1\}}$, donn\'e par 
l'\'equation $w^p-\bar u$, et $d\bar u$ ne 
poss\`ede ni z\'ero, ni p\^ole sur cet ouvert, donc 
$\spec \bar \mathcal B$ est lisse, et c'est un ouvert de la droite projective 
sur $k$ de la forme $\mathbb P_k^1 \setminus \{\infty,z_a\}_{a \in \aret (r_1)\setminus \{a_1\}}$, 
o\`u les $z_a$ sont des 
\'el\'ements de $k$ distincts deux \`a deux. Autrement dit, on a  
$\bar \mathcal B=k[z,(z-z_a)^{-1}]_{a \in \aret (r_1)\setminus \{a_1\}}$. 
Soit $Z$ (resp. $Z_{a}$) un rel\`evement de $z$ (resp. $z_{a}$) dans 
$\mathcal B$ (resp. $R$), comme $\mathcal B$ est plat sur $R$, on a 
$\mathcal B=R\{Z,(Z-Z_{a})^{-1}\}_{a \in \aret (r_1)\setminus \{a_1\}}$.
\par
Pour $a \in \aret (r_1)\setminus \{a_1\}$, on a 
$1+\ord_{t=t_a}d\bar u=-m(a)\ne 0 \mod p$ ; il r\'esulte 
alors de la proposition \ref{reducbord} (B) que 
$\mathcal B \otimes_{\mathcal A}R[[T-T_{a}]]\{(T-T_{a})^{-1})\}$ 
est \'egal \`a $R[[Z-Z_{a}]]\{(Z-Z_{a})^{-1}\}$, et l'action induite par 
celle de $\sigma$ est donn\'ee sur un param\`etre convenable 
$\tilde X_a$ de $R[[Z-Z_{a}]]\{(Z-Z_{a})^{-1}\}$ par 
$$\sigma (\tilde X_a)=\zeta^{-\frac 1m(a)}\tilde X_a
(1+\pi^{v_K(\lambda)-n_1}\tilde X_a^{m(a)})^{-\frac 1m(a)}.$$
\par
Soit $\psi_a$ l'isomorphisme \'equivariant de 
$\mathcal O(\mathcal D_a)_{(\pi)}^\wedge=R[[X_a]]\{X_a^{-1}\}$ sur \\
$R[[Z-Z_{a}]]\{(Z-Z_{a})^{-1}\}$ qui envoie $X_a$ sur $\tilde X_a$. 
En utilisant le lemme \ref{patchingcloseddisc}, on construit 
alors une $R$-alg\`ebre de la forme $R\{Z_0\}$, muni d'un $R$-automorphisme 
$\sigma$ d'ordre $p$, dont 
l'arbre de Hurwitz associ\'e s'identifie \`a $(\Gamma, \mathcal H)$. De plus, 
comme $\ord_{\infty}d\bar u=-(m(a_1)+1)$, la proposition \ref{reducbord} (B) 
montre l'existence d'un param\`etre $V$ v\'erifiant les conditions 
attendues.
\end{proof}
Comme pr\'ec\'edemment, en en d\'eduit le
\begin{cor}\label{hoadd}
Soit $\Gamma$ un arbre fini connexe et $\mathcal H=(r_0,d_0,\epsilon,m,h)$ 
une donn\'ee de Hurwitz sur $\Gamma$ ; on suppose que $(\Gamma, \mathcal H)$ 
v\'erifie les conditions suivantes :
\begin{itemize}
\item
La racine $r_0$ est additive ou \'etale.
\item
Il existe une unique ar\^ete $a_0$ d'origine $r_0$.
\item
Si $a$ est une ar\^ete d'origine $r_1:=t_\Gamma(a_0)$ distincte de 
$a_1:=\bar a_0$, l'arbre de 
Hurwitz $(\Gamma_a,\mathcal H_a)$ 
provient d'un $R$-automorphisme $\sigma_a$ d'ordre $p$ d'un disque formel, 
tel qu'il existe un param\`etre $X_a$ avec 
$$\sigma_{a}(X_a)=\zeta^{-\frac 1{m(a)}}X_a
(1+\pi^{v(\lambda)-n_1}X_a^{m(a)})^{-\frac 1{m(a)}},$$
o\`u 
$n_1$ est l'entier d\'efini par la relation $v_K(p)=d(r_1)+n_1(p-1)$.
\item
Le sommet $r_1$ est r\'ealisable.
\end{itemize}
Alors, si $R$ d\'esigne l'anneau de valuation de $K$, 
il existe un automorphisme $\sigma$ d'ordre $p$ du disque 
formel $\spec R[[Z]]$ dont l'arbre de Hurwitz associ\'e est isomorphe \`a 
$(\Gamma, \mathcal H)$, et tel qu'il existe un param\`etre $X$ au bord du 
disque formel avec 
$$\sigma(X)=\zeta^{-\frac 1{m(a_0)}} 
X(1 +\pi^{v_K(\lambda)-n_0}X^{m(a_o)})^{-\frac 1{m(a_0)}},$$
o\`u $n_0$ est l'entier d\'efini par la relation $v_K(p)=d_0+(p-1)n_0$.
\end{cor}
\subsection[D\'emonstration]{D\'emonstration du th\'eor\`eme \ref{hurwitzdisc}}
On raisonne par r\'ecurrence sur le nombre maximal $N\ge 3$ de sommets d'une 
cha\^ \i ne d'origine $r_0$ et de sommet terminal maximal.
Notons $P_N$ la propri\'et\'e suivante : Tout arbre de Hurwitz v\'erifiant 
les conditions $D[i]$, o\`u $1 \le i \le 3$, et dont le nombre 
maximum de sommets d'une 
cha\^ \i ne d'origine $r_0$ et de sommet terminal maximal est inf\'erieur ou 
\'egal \`a $N$, provient d'un automorphisme $\sigma$ du disque formel, tel 
qu'il existe un param\`etre $X$ au bord du disque, avec 
$\sigma(X)=\zeta^{-\frac 1{m(a_0)}} 
X(1 +\pi^{v_K(\lambda)-n_0}X^{m(a_o)})^{-\frac 1{m(a_0)}}$.
\par 
La propri\'et\'e $P_3$ r\'esulte du corollaire \ref{homult}.
\par
Si maintenant $N \ge 4$, on suppose $P_{N-1}$ d\'emontr\'ee ; soit 
$(\Gamma, \mathcal H)$ un arbre de Hurwitz v\'erifiant les conditions 
$D[i]$, o\`u $1 \le i \le 3$, et tel que le nombre maximal 
de sommets d'une 
cha\^ \i ne d'origine $r_0$ et de sommet terminal maximal est \'egal \`a $N$. 
Alors, les hypoth\`eses du corollaire \ref{hoadd} sont 
v\'erifi\'ees, et donc l'arbre de Hurwitz $(\Gamma, \mathcal H)$ provient 
d'un automorphisme d'ordre $p$ d'un disque formel, avec un bon param\`etre au 
bord.
\subsection{Un crit\`ere de r\'ealisabilit\'e}
\subsubsection{Notations}
\par
Soit $I$ un ensemble fini, non vide, on note $\mathcal P_I$ l'ensemble des 
partitions de $I$. L'ensemble $\mathcal P_I$ est en bijection avec l'ensemble 
des relations d'\'equivalence sur $I$ : Plus pr\'ecis\'ement, si 
$\mathcal R$ est une relation d'\'equivalence sur $I$, 
la partition associ\'ee est par d\'efinition l'ensemble des classes 
d'\'equivalence ; r\'eciproquement, si $P$ est une partition, on d\'efinit 
une relation d'\'equivalence sur $I$ par $i\equiv j \mod P$ si et seulement si 
il existe $J$ dans $P$ tel que $i$ et $j$ appartiennent \`a $J$. La partition 
associ\'ee \`a la relation d'\'egalit\'e sur $I$ sera not\'ee $P_=$.
\par
Il existe une relation d'ordre naturelle sur $\mathcal P_I$, \`a savoir 
$P<P'$ si et seulement si pour tout $J'$ dans $P'$, il existe $J$ dans $P$ 
tel que $J' \subset J$. Autrement dit, si la relation d'\'equivalence 
correspondant \`a $P'$ est plus fine que celle correspondant \`a $P$.
\subsubsection{Le $k$-sch\'ema $X_{\underline e}$.}
Soient $m$ un entier positif premier \`a $p$, $I$ un ensemble fini \`a 
$m+1$ \'el\'ements, 
$\underline e:=(e_i)_{i \in I}$ une famille d'\'el\'ements non nuls de 
$\mathbb F_p$, telle que $\sum_{i \in I}e_i=0$ ; on notera 
$X_{\underline e}$ le sous-sch\'ema ferm\'e de $\spec k[T_i]_{i \in I}$ 
d'id\'eal $(\sum_{i \in I}e_iT_i^\nu)_{1 \le \nu \le m-1, (\nu,p)=1}$. 
Soit $x:=(t_i)$ un point rationnel sur $k$ de $X_{\underline e}$, il 
d\'efinit une partition $P_x$ de $I$ telle que $i=j \mod P_x$ si et 
seulement si $t_i=t_j$. Pour une partition $P$ fix\'ee de $I$, on notera 
$X_{\underline e,P}$ le sous-sch\'ema ferm\'e de $X_{\underline e}$ d'id\'eal 
$(T_i-T_j)_{i\equiv j \mod P}$. En particulier, 
$X_{\underline e,P_=}=X_{\underline e}$. On notera $X_{\underline e,P}^*$ 
l'ouvert de $X_{\underline e,P}$ compl\'ementaire de la r\'eunion des 
$X_{\underline e,P'}$, avec $P'<P$. On remarquera que 
$X_{\underline e,P}^*(k)=\{x \in X_{\underline e}(k), P_x=P\}$
\begin{lem}\label{vanderMonde}
Soit $P$ une partition de $I$, diff\'erente de $P_=$. alors si 
$x=(t_i)$ est un point rationnel sur $k$ de $X_{\underline e}$ avec $P_x=P$, 
on a $\sum_{i \in J}e_i=0$ pour tout $J$ dans $P$.
\end{lem}
\begin{proof}
\par
Pour $J$ dans $P$, notons $a_J=t_i$ pour n'importe quel $i$ dans $J$, ce qui 
est bien d\'efini car $P=P_x$. Alors, on a 
$\sum_{J \in P}(\sum_{i \in J}e_i)a_J^{\nu}=0$ pour $\nu$ variant de $0$ \`a 
$m-1$. Comme  $P$ est diff\'erente de $P_=$, le cardinal $r$ de $P$ est 
inf\'erieur ou \'egal \`a $m$. Les $a_J$ \'etant distincts deux \`a deux, 
on d\'eduit des $r$ premi\`eres \'equations la relation 
$\sum_{i \in J}e_i=0$ pour $J$ dans $P$.
\end{proof}
\begin{defn}
On dira qu'une partition de $I$ est $\underline e$-adapt\'ee si pour tout 
$J$ dans $P$, on a $\sum_{i \in J}e_i=0 \mod p$. On notera 
$\mathcal P_{I,\underline e}$ l'ensemble des partitions 
$\underline e$-adapt\'ees de $I$.
\end{defn}
\begin{cor}\ \\
(a) Si $P$ est une partition de $I$, $\underline e$-adapt\'ee, alors 
$X_{\underline e,P}$ est isomorphe \`a l'espace affine 
$\mathbb A_k^P:=\spec k[A_J]_{J \in P}$. En particulier, 
$X_{\underline e,P}$ est de dimension $\card P$.\\
(b) Soit $Y_{\underline e}:=
\cup_{P \in \mathcal P_{I,\underline e}}X_{\underline e,P}$, 
alors $X_{\underline e, \text{r\'ed}}$ est la r\'eunion disjointe de 
$Y_{\underline e}$ et de \text{$X^*_{\underline e}:=X_{\underline e,P_=}^*$.}
\end{cor}
\begin{proof}\ \\
(a) L'isomorphisme en question est donn\'e par l'isomorphisme
$$\frac{k[T_i]_{i \in I}}{(T_i-T_j)_{i\equiv j \mod P}} \to k[A_J]_{J \in P}$$
qui envoie $T_i$ sur $A_J$, o\`u $J={i \mod P}$.\\
(b) r\'esulte du lemme \ref{vanderMonde}, qui entra\^ \i ne que 
$X_{\underline e, \text{r\'ed}}(k)$ est la r\'eunion disjointe de 
$Y_{\underline e}(k)$ et de $X^*_{\underline e}(k)$.
\end{proof}
\begin{prop}\ \\
\label{criteresuffisant}
(a) Les composantes irr\'eductibles de $Y_{\underline e}$ sont les 
$X_{\underline e,P}$, o\`u $P$ est une partition maximale 
$\underline e$-adapt\'ee.\\
(b) Si il existe une partition maximale $\underline e$-adapt\'ee avec 
$\card P \le [\displaystyle{\frac mp}]+1$, alors $X^*_{\underline e}(k)$ 
est non vide.\\
(c) Si $X^*_{\underline e}(k)$ est non vide, l'ouvert $X^*_{\underline e}$ 
est r\'egulier.
\end{prop}
\begin{proof}
\par
L'assertion (a) est claire. Pour toute composante irr\'eductible $Z$ de 
$X_{\underline e}$, la dimension de $Z$ est sup\'erieure ou \'egale \`a 
$m+1-(m-1-[\displaystyle{\frac mp}])=[\displaystyle{\frac mp}]+2$. Soit 
$P$ une partition maximale $\underline e$-adapt\'ee avec 
$\card  P \le [\displaystyle{\frac mp}]+1$, alors $X_{\underline e,P}$ 
n'est pas une composante irr\'eductible de $X_{\underline e}$, et une 
composante irr\'eductible de $X_{\underline e}$ contenant 
$X_{\underline e,P}$ doit alors rencontrer $X^*_{\underline e}(k)$. Il reste 
\`a voir l'assertion (c) : On voit imm\'ediatement que l'espace tangent 
en un point ferm\'e de $X^*_{\underline e}$ est de dimension 
$[\displaystyle{\frac mp}]+2$, donc inf\'erieur \`a la dimension en ce point. 
On applique alors le crit\`ere jacobien.
\end{proof}
\begin{rem}
\label{dimension}
Le stabilisateur $S_{\infty}:=\{z \mapsto az+b|a \in k^*, b \in k\}$ 
dans PGL$(2,k)$ du point $\infty$ dans $\mathbb P^1(k)$ est de dimension 2, 
et agit librement sur $X^*_{\underline e}$. Si 
$X^*_{\underline e}(k)$ est non vide, le quotient de 
$X^*_{\underline e}$ par l'action de $S_{\infty}$ est donc une vari\'et\'e 
alg\'ebrique affine sur $k$ de dimension $[\displaystyle{\frac mp}]$.
\end{rem}
\subsubsection{Crit\`ere de r\'ealisabilit\'e}
\begin{prop}\label{critrealdisc}
\par
Soit $(\Gamma,\mathcal H)$ un arbre de Hurwitz, et $s$ un sommet 
de $\Gamma$.\\
(a) Supposons que $s$ est multiplicatif, et qu'il existe une 
unique ar\^ete $a_0$ d'origine $s$ telle que $m(a_0) \ne 0$. 
Notons $m:=-m(a_0)$, $I_s:=\aret(s) \setminus a_0$ et 
$\underline h:=(h(a))_{a \in I_s}$.
\par
Alors, si il existe une 
partition maximale $\underline h$-adapt\'ee $P$ de $I_s$ avec\\
$\card P\le [\displaystyle{\frac {m}p}]+1$, le 
sommet $s$ est r\'ealisable.\\
(b) Supposons que $s$ est additif, et qu'il existe une 
unique ar\^ete $a_0$ d'origine $s$ telle que $m(a_0)<0$. 
Notons $m:=-m(a_0)$, $I_s$ la r\'eunion disjointe de 
$I'_s:=\aret(s) \setminus a_0$ 
et d'un ensemble fini \`a 
$m+1-\card(I'_s)$ \'el\'ements, et 
$\underline e=(e_i)_{i \in I_s}$ la famille d\'efinie par 
$e_i:=m(i)$ si $i \in I'_{s}$ et $e_i=-1$ sinon.
\par
Alors, si il existe une 
partition maximale $\underline e$-adapt\'ee $P$ de $I_s$ avec\\
$\card P\le [\displaystyle{\frac {m}p}]+1$, le 
sommet $s$ est r\'ealisable.
\end{prop}
\begin{rem}
La proposition pr\'ec\'edente g\'en\'eralise la proposition III 4.5.1 de 
\cite{G-M 2}. Plus pr\'ecis\'ement, avec la terminologie employ\'ee dans 
cette proposition, le fait que $(h_i)$ appartienne au lieu critique revient 
\`a dire que la partition grossi\`ere est maximale, ce qui pour $m<p$ est 
\'equivalent au crit\`ere \ref{critrealdisc}. Par ailleurs, l'exemple 
III 4.7 de ce m\^eme article prouve que le crit\`ere n'est pas n\'ecessaire.
\end{rem}
\begin{proof}\ \\
(a) On cherche $\bar u:=\bar u_s$ sous la forme 
$\bar u=\displaystyle{\prod_{i \in I_s}(t-t_i)^{h_i}}$, 
en convenant que \text{$j_s(a_0)=\infty$} et que $h_i$ d\'esigne 
un repr\'esentant de $h(i)$ dans $\mathbb Z$, avec la relation 
$\sum_{i \in I_s}h_i=0$. Soit $x:=t^{-1}$, on a alors 
$\bar u=\prod_{i \in I_s}(1-t_ix)^{h_i}$, d'o\'u 
$$\frac {d\bar u}{\bar u}=\sum_{\nu=0}^{+\infty}
(\sum_{i \in I_s}h_it_i^{\nu+1})x^{\nu}dx.$$
La relation $\Div(\displaystyle{\frac{d\bar u_s}{\bar u_s}})=-\sum_{a \in \aret(s)}(m(a)+1)[j_s(a)]$ est alors \'equivalente \`a 
$\ord_\infty\displaystyle{\frac {d\bar u}{\bar u}}=m-1$, ce qui 
revient \`a dire que $(t_i)$ est un point rationnel sur $k$ de 
$X^*_{\underline h}$. On utilise alors la proposition 
\ref{criteresuffisant}.\\
(b) On cherche $\bar u:=\bar u_s$ sous la forme 
$$\bar u=\frac{\prod_{i\in I_s \setminus I'_s}(t-t_i)}
{\prod_{i \in I'_s}(t-t_i)^{m(i)}}.$$
On convient que $j_s(a_0)=\infty$. Comme ci-dessus, on montre 
alors que la relation
$$\Div(d\bar u_s)=-\sum_{a \in \aret(s)}(m(a)+1)[j_s(a)]$$
revient \`a dire que $(t_i)$ est un point rationnel sur $k$ de 
$X^*_{\underline e}$.
\end{proof}
\begin{cor}
Si $p=2$, tout sommet non maximal d'un arbre de Hurwitz 
est r\'ealisable.
\end{cor}
\par
Ce dernier r\'esultat peut se d\'emontrer aussi directement, en 
exhibant un point de $X_{\underline e}^*(k)$, avec $I$ un 
ensemble de cardinal pair et $e_i=1$ pour 
tout $i$ dans $I$.
\subsection{Un exemple avec $p=5$}
\par
On donne ici un exemple d'arbre de Hurwitz provenant d'un automorphisme 
d'ordre $5$ d'un disque formel sur un anneau de valuation discr\`ete complet 
d'in\'egales caract\'eristiques $(0,5)$, sans toutefois pr\'eciser 
la m\'etrique. Par souci de concision, sur la figure n'apparaissent 
que les valeurs non nulles de $h$ et de $m$. 
Les ar\^etes qui ne sont pas des feuilles sont indiqu\'ees en 
gras sur la figure.

\begin{figure}[htpb]
         \input{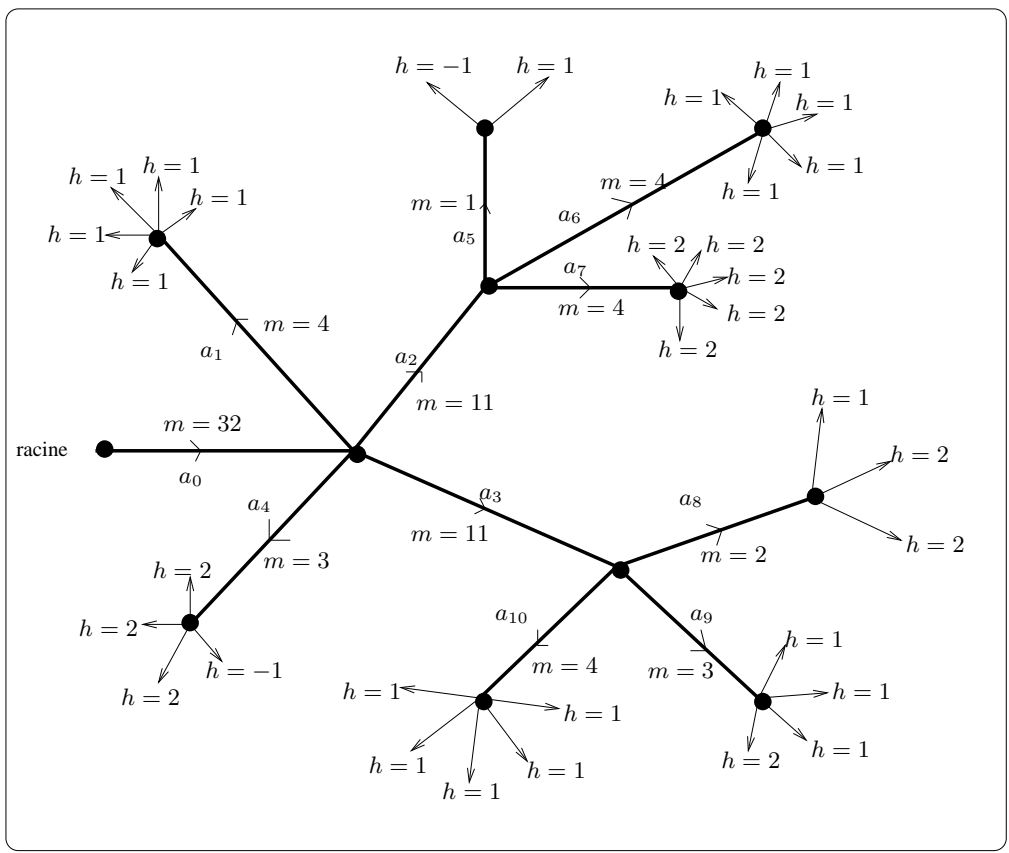}
         \caption{Un arbre de Hurwitz r\'ealisable pour $p=5$}
\end{figure}
\begin{prop}
Soit $K_0$ le corps des fractions de $W(k)[\zeta]$, o\`u $k$ est un 
corps alg\'ebriquement clos de caract\'eristique $5$ et $\zeta$ une racine 
primitive $5$-i\`eme de l'unit\'e. 
Il existe une extension finie $K$ de $K_0$ et une m\'etrique $\epsilon$ sur 
l'arbre ci-dessus telle que l'arbre de Hurwitz obtenu provienne 
d'un automorphisme d'ordre $5$ du disque formel sur l'anneau de valuation $R$ 
de $K$, op\'erant sans inertie au bord, avec conducteur de Hasse $32$.
\end{prop}
\begin{proof}
\par
Posons $K:=K_0[\lambda^{\frac 1{4224}}]$, o\`u $\lambda=\zeta-1$. Dans ce cas, 
on peut par exemple prendre la m\'etrique suivante (bien \'evidemment, c'est 
tr\`es loin d'\^etre la seule possibilit\'e pour une extension de $K_0$ aussi ramifi\'ee) : 
Les \'epaisseurs des ar\^etes $a_i$ sont : 
$\epsilon(a_0):=33$, 
$\epsilon(a_1):=792$, 
$\epsilon(a_2):=96$, 
$\epsilon(a_3):=96$, 
$\epsilon(a_4):=1056$, 
$\epsilon(a_5):=2112$, 
$\epsilon(a_6):=528$, 
$\epsilon(a_7):=528$, 
$\epsilon(a_8):=1056$, 
$\epsilon(a_9):=704$, 
$\epsilon(a_{10}):=528$.
\par
Avec un tel $K$ et la m\'etrique d\'efinie ci-dessus (avec la racine \'etale), on obtient un arbre de 
Hurwitz en utilisant la d\'efinition de $m$ et de $h$ donn\'ee sur le dessin. 
Il reste \`a voir que les sommets de valence sup\'erieure \`a $3$ sont 
r\'ealisables. Pour cela, on applique la proposition \ref{critrealdisc} :\\
Les sommets multiplicatifs sont tous r\'ealisables, car la partition 
grossi\`ere est maximale dans chaque cas.\\
Le sommet additif origine de $a_5$ est r\'ealisable : La partition 
$\underline e$-adapt\'ee
$$\left\{ (4,-1,-1,-1,-1),(4,-1,-1,-1,-1),(1,-1)\right\}$$
de $\underline e:=(4,4,1,-1,-1,-1,-1,-1,-1,-1,-1,-1)$ est maximale, 
de cardinal $3=[\frac{11}5]+1$.\\
Le sommet additif origine de $a_8$ est r\'ealisable : La partition 
$\underline e$-adapt\'ee
$$\left\{ (4,-1,-1,-1,-1),(3,-1,-1,-1),(2,-1,-1)\right\}$$
de $\underline e:=(2,3,4,-1,-1,-1,-1,-1,-1,-1,-1,-1)$ est maximale, de cardinal 
$3=[\frac{11}5]+1$.\\
Le sommet additif origine de $a_2$ est r\'ealisable : La partition 
$\underline e$-adapt\'ee
$$\left\{ (1,1,3),(4,-1,-1,-1,-1),(-1,-1,-1,-1,-1),\dots,
(-1,-1,-1,-1,-1)\right\}$$
(o\`u $(-1,-1,-1,-1,-1)$ est r\'ep\'et\'e $5$ fois) de 
$\underline e:=(1,1,3,4,4,-1,\dots,-1)$ (o\`u $-1$ est r\'ep\'et\'e 
$29$ fois) est maximale, de cardinal 
$7=[\frac{32}5]+1$.\\
\end{proof}
\section{Automorphisme d'ordre $p$ d'une couronne formelle}

\subsection{Mod\`ele minimal d'une couronne ouverte avec action d'un automorphisme d'ordre $p$}

\par
Consid\'erons un $R$-automorphisme $\sigma$ d'ordre $p$ de la 
couronne formelle $\mathcal C_e$, 
$F_{\sigma}$ l'ensemble de ses points fixes g\'eom\'etriques dans 
$\mathcal C_{e,K}:=\mathcal C_e \times_{\spec R}\spec K$, (\'eventuellement 
vide). De plus, on supposera tous les points fixes rationnels sur $K$. Il 
existe un unique mod\`ele minimal $\mathcal M(\mathcal C_{e,K},\sigma)$ 
de $\mathcal C_{e,K}$ 
pour lequel les sp\'ecialisations des points de $F_{\sigma}$ sont lisses 
et distinctes : Si $F_{\sigma}$ est vide, 
$\mathcal M(\mathcal C_{e,K},\sigma)$ est \'egal \`a $\mathcal C_e$. Sinon, 
il est construit de la fa\c con suivante : Choisissons une 
coordonn\'ee de Laurent $Z$ sur $\mathcal C_e$, on commence par 
\'eclater les id\'eaux $(\pi^{a},Z)$, o\`u 
$a$ d\'ecrit l'ensemble des valuations prises par les points de $F_{\sigma}$. 
On obtient ainsi une 
cha\^ \i ne de droites projectives reli\'ee \`a chaque extr\'emit\'e \`a la 
transform\'ee stricte d'un des deux bords de $\mathcal C_e$. Si le mod\`ele 
ainsi obtenu ne s\'epare pas les 
sp\'ecialisations des points de $F_{\sigma}$, on \'eclate encore les fibres 
formelles contenant plus d'un point de $F_{\sigma}$ (ce sont des disques 
formels) pour obtenir le mod\`ele minimal cherch\'e. Sa fibre sp\'eciale 
est un arbre de droites projectives, reli\'e aux composantes correspondant aux 
transform\'es strictes des bords. La cha\^ \i ne de 
composantes irr\'eductibles obtenue dans la premi\`ere \'etape, 
menant d'un bord \`a l'autre, s'appelle la 
{\bf  cha\^ \i ne fondamentale}. Les composantes de la 
cha\^ \i ne fondamentale seront dites fondamentales. 
\subsection{Arbre de Hurwitz associ\'e \`a un automorphisme d'ordre $p$ d'une 
couronne formelle}
\par
Soit $\Gamma^*_{\sigma}$ l'arbre dual de la fibre sp\'eciale de 
$\mathcal M_{\sigma}:=\mathcal M(\mathcal C_{e,K},\sigma)$. Choisissant 
l'un des deux bords 
$\eta$ de $\mathcal C_e$, on oriente $\Gamma^*_{\sigma}$ \`a partir du sommet 
$r_\eta$ correspondant. Soit $\Gamma_{\sigma}$ l'arbre orient\'e 
d\'ecrit ci-dessous :
\begin{itemize}
\item
L'ensemble des sommets de $\Gamma_{\sigma}$ est la r\'eunion disjointe de 
$\som \Gamma^*_{\sigma}$ et de $F_{\sigma}$. Le sommet correspondant \`a un 
point $x$ de $F_{\sigma}$ sera not\'e $s_x$.
\item
L'ensemble des ar\^etes positives de $\Gamma_{\sigma}$ est la r\'eunion 
disjointe de 
$\aret^+ \Gamma^*_{\sigma}$ et de $F_{\sigma}$. L'ar\^ete positive 
correspondant \`a un point $x$ de $F_{\sigma}$ sera not\'e $a_x$.
\item
Si $x$ est un point de $F_{\sigma}$, le sommet terminal de $a_x$ est $s_x$ 
et l'origine de $a_x$ est le sommet de $\Gamma^*_{\sigma}$ sur lequel se 
sp\'ecialise $x$.
\end{itemize}
\par
Soit $s$ un sommet de $\Gamma^*_{\sigma}$ ; le point g\'en\'erique de 
la composante irr\'eductible 
correspondant \`a $s$ d\'efinit une valuation discr\`ete $v_s$ du corps des 
fractions $\mathcal K_e$ de $\mathcal A_e$. On note $\delta(s)$ la valuation 
de la 
diff\'erente de l'extension de corps valu\'es $(\mathcal K_e,v_s)/
(\mathcal K_e^{\sigma},v_s)$. D'apr\`es la proposition \ref{reducmuptorseur}, 
on a $0 \le \delta(s) \le v_K(p)$. Un sommet $s$ de $\Gamma^*_{\sigma}$ sera 
dit multiplicatif si $\delta(s)=v_K(p)$, additif si\\
$0<\delta(s)<v_K(p)$ et \'etale si $\delta(s)=0$.
\par
On d\'efinit une application $\epsilon:=\epsilon_{\sigma}$ de 
$\aret \Gamma_{\sigma}$ dans $\mathbb N$ : Si $a$ est une ar\^ete de 
$\Gamma^*_{\sigma}$, elle correspond \`a un point double $z_a$ de 
$\mathcal M(\mathcal C_{e,K},\sigma)$, on d\'efinit alors $\epsilon(a)$ 
comme \'etant l'\'epaisseur de $z_a$ dans 
$\mathcal M(\mathcal C_{e,K},\sigma)$. Si maintenant $a=a_x$, o\`u $x$ est 
un point de $F_{\sigma}$), on d\'efinit $\epsilon(a):=0$.
\par
On d\'efinit une application $m:=m_{\sigma}$ de $\aret \Gamma_{\sigma}$ dans 
$\mathbb Z$ et une application $h:=h_{\sigma}$ de $\aret \Gamma_{\sigma}$ dans 
$\mathbb Z/p\mathbb Z$ de la fa\c con suivante : Soit $a$ une ar\^ete de 
$\Gamma_{\sigma}$. 
Si $a$ est une ar\^ete de 
$\Gamma^*_{\sigma}$ (resp. $a=a_x$, o\`u $x$ est un point de $F_{\sigma}$), 
elle correspond \`a un point double orient\'e $z_a$ de 
$\mathcal M(\mathcal C_{e,K},\sigma)$ (resp. \`a la sp\'ecialisation de $x$ 
dans la fibre sp\'eciale de $\mathcal M(\mathcal C_{e,K},\sigma)$). 
Soit $\xi_a$ le bord de la couronne formelle (resp. du disque formel) 
$\spec \hat {\mathcal O_{\mathcal M_{\sigma},z_a}}$ qui correspond \`a 
l'origine de $a$. On notera $\hat {\mathcal O_a}$ 
le localis\'e-compl\'et\'e 
$(\hat {\mathcal O_{\mathcal M_{\sigma},z_a})}_{\xi_a}^{\wedge}$. On note 
$\omega_a$ la $1$-forme diff\'erentielle associ\'ee au torseur 
$\spec \hat {\mathcal O_a} \to \spec \hat {\mathcal O_a^{\sigma}}$ 
si ce torseur est radiciel en r\'eduction. 
On pose $m(a):=-(1+\ord_{\xi_a}\omega_a)$ et 
$h(a)$ \'egal au r\'esidu de $\omega_a$ si l'origine de $a$ est un sommet 
multiplicatif. Si maintenant l'origine de $a$ est 
un sommet additif, on pose 
$m(a)=-(1+\ord_{\xi_a}\omega_a)$ et $h(a)=0$. Si l'origine de $a$ est un 
sommet \'etale, $m(a)$ est le conducteur de Hasse de l'extension 
$\hat {\mathcal O_a}\otimes_R k/\hat {\mathcal O_a}^{\sigma}\otimes_R k$. 
Enfin, $m(\bar a_x)=-m(a_x)$ 
et $h(\bar a_x)=-h(a_x)$.
\begin{prop}\ \\
\label{couronnehurwitzdata}
\par
La donn\'ee 
$\mathcal H_{\sigma}:=(r_{\eta},\delta(r_{\eta}),\epsilon_{\sigma},
m_{\sigma}, h_{\sigma})$ est une donn\'ee de Hurwitz sur 
$\Gamma_{\sigma}$. On dira que l'arbre de Hurwitz 
$(\Gamma_{\sigma},\mathcal H_{\sigma})$ est l'arbre de Hurwitz associ\'e 
\`a l'automorphisme $\sigma$ de 
la couronne formelle d'\'epaisseur $e$. Par ailleurs,
\begin{itemize}
\item
(Loi de variation de la diff\'erente) Pour tout sommet $s$ de 
$\Gamma_\sigma$, on a
$$\delta(s)=d(s):=
\delta(r_{\eta})+(p-1)\sum_{a \in \aret^+\Gamma_s}m(a)\epsilon(a),$$
o\`u $\Gamma_s$ est la cha\^ \i ne d'origine $r_{\eta}$ et de sommet 
terminal $s$.
\item
Le cardinal de $F_{\sigma}$ est \'egal \`a $m(a_{\eta})+m(a_{\eta'})$, o\`u 
$a_{\eta}$ (resp. $a_{\eta'}$) est l'unique ar\^ete d'origine $r_{\eta}$ 
(resp. $r_{\eta'}$).
\item
Les sommets de valence sup\'erieure ou \'egale \`a $3$ sont r\'ealisables.
\item
Un sommet maximal est soit le sommet terminal 
d'une feuille, soit $r_{\eta'}$, o\`u $\eta'$ est le bord de la couronne 
oppos\'e \`a $\eta$ ($r_{\eta'}$ est le sommet terminal de la cha\^ \i ne 
fondamentale).
\item
Si $a$ est une feuille de $(\Gamma_{\sigma},\mathcal H_{\sigma})$, alors 
l'origine de $a$ est soit un sommet maximal de $\Gamma_{\sigma}^*$, soit un 
sommet fondamental.
\item
Il existe au plus deux sommets fondamentaux qui sont multiplicatifs. De plus, 
s'il y en a deux, il existe une ar\^ete qui les relie.
\end{itemize}
\end{prop}
\begin{proof}
Le fait que $\mathcal H_{\sigma}$ soit une donn\'ee de Hurwitz se montre 
exactement de la m\^eme fa\c con que pour le disque. En particulier, on 
prouve de m\^eme que les sommets de valence sup\'erieure ou \'egale 
\`a $3$ sont r\'ealisables. La loi de la variation de la diff\'erente 
r\'esulte de la proposition \ref{pentediff}.
\par
Soit $s$ un sommet maximal. Supposons $s$ non fondamental, 
on note alors $a$ l'unique ar\^ete de sommet origine 
fondamental et de sommet terminal non fondamental qui apparait dans la 
cha\^ \i ne positive $\Gamma_s$. L'arbre de Hurwitz 
$(\Gamma_{\sigma}[a],\mathcal H_{\sigma}[a])$ provient alors d'un 
automorphisme d'ordre $p$ du disque formel, on sait alors que $s$ est le 
sommet terminal d'une feuille de 
$(\Gamma_{\sigma}[a],\mathcal H_{\sigma}[a])$, donc de 
$(\Gamma_{\sigma},\mathcal H_{\sigma})$. Si maintenant $s$ est fondamental, 
comme il est maximal il doit \^etre \'egal \`a $r_{\eta'}$.
\par
Montrons \`a pr\'esent la formule donnant le cardinal de $F_{\sigma}$, 
c'est-\`a-dire le nombre de feuilles de l'arbre de Hurwitz 
$(\Gamma_{\sigma},\mathcal H_{\sigma})$. Cela peut se voir par r\'ecurrence 
sur le nombre de sommets fondamentaux :\\
Si il n'y a que deux sommets fondamentaux, $F_{\sigma}$ est vide, et la 
proposition \ref{pentediff} montre la formule.\\
Soit $N \ge 3$, supposons la formule d\'emontr\'ee pour les 
automorphismes d'ordre $p$ de couronnes 
formelles dont l'arbre de Hurwitz poss\`ede $N-1$ sommets fondamentaux. 
Soit $\sigma$ un automorphisme d'ordre $p$ d'une couronne 
formelle $\mathcal C$ d'\'epaisseur $e$ dont l'arbre de Hurwitz poss\`ede 
$N$ sommets fondamentaux. 
Alors, si $s$ est le sommet fondamental sommet terminal de l'ar\^ete 
$a_\eta$ et $a$ l'unique ar\^ete fondamentale d'origine $s$, 
l'arbre de Hurwitz $(\Gamma_\sigma[a], \mathcal H_\sigma[a])$ provient de 
l'automorphisme $\sigma$ restreint \`a la couronne formelle 
$\mathcal C[a]$ de fibre g\'en\'erique form\'ee des points $x$ tels que 
$\epsilon(a)<v_K(Z(x))<e$, o\`u $Z$ est 
une coordonn\'ee de Laurent sur $\mathcal C$. Le nombre de feuilles de 
$(\Gamma_\sigma[a], \mathcal H_\sigma[a])$ est alors $m(a)+m(a_{\eta'})$. 
Si $b$ est une ar\^ete non fondamentale de $\Gamma_\sigma$ d'origine $s$, 
le nombre de feuilles de $(\Gamma_\sigma[b], \mathcal H_\sigma[b])$, qui 
provient d'un automorphisme de disque formel, est \'egal \`a $m(b)+1$. 
Donc, si $I$ d\'esigne l'ensemble des ar\^etes non fondamentales d'origine 
$s$, le cardinal de $F_{\sigma}$ est, d'apr\`es $H[3]$ :
\begin{eqnarray*}
\sum_{b \in I}(m(b)+1)+m(a)+m(a_{\eta'})
= & &-(m(\bar a_\eta)+1)+ 2-1+m(a_{\eta'})\\
= & &m(a_\eta)+m(a_{\eta'})
\end{eqnarray*}
\par
Il reste \`a voir la derni\`ere assertion. 
Pour $s$ un sommet fondamental distinct de $r_{\eta}$ et $r_{\eta'}$, soit 
$a_1$ l'ar\^ete fondamentale de sommet terminal $s$ et $a_2$ l'ar\^ete 
fondamentale d'origine $s$. Alors la pente de la variation de la diff\'erente 
le long de la cha\^ \i ne fondamentale est \'egale \`a $m(a_1)(p-1)$ 
avant le sommet $s$ et \`a $m(a_2)(p-1)$ apr\`es. Donc, la variation de 
la pente au passage de $s$ est \'egale, d'apr\`es $H[3]$, \`a 
$-\sum_{a \in \aret(s), a\ne a_1,a_2}(m(a)+1)<0$, car pour une ar\^ete $a$ non 
fondamentale d'origine $s$, $(\Gamma_{\sigma}[a],\mathcal H_{\sigma}[a])$ 
provient d'un automorphisme d'ordre $p$ du disque formel et donc $m(a) \ge 0$. 
La variation de la diff\'erente le long de la cha\^ \i ne fondamentale est 
donc une fonction continue, lin\'eaire par morceaux et concave (i.e. la pente 
diminue). Supposons qu'il existe un sommet fondamental multiplicatif, et 
prenons $s$ le plus petit d'entre eux. Alors, soit $s$ est maximal, et alors 
c'est le seul sommet multiplicatif fondamental, soit il existe une (unique) 
ar\^ete fondamentale $a$ d'origine $s$ et alors $m(a)\ge 0$. Si $m(a)=0$, le 
sommet terminal de $a$ est multiplicatif, et ensuite la diff\'erente doit 
chuter, donc les \'eventuels sommets suivants ne sont pas multiplicatifs. 
Sinon, $s$ est le seul sommet multiplicatif fondamental.
\end{proof}
\subsection{Th\'eor\`eme de r\'ealisation}
\par
Comme pour le cas du disque, on peut donner une caract\'erisation des 
arbres de Hurwitz provenant d'un automorphisme d'ordre $p$ d'une 
couronne formelle (ne permutant pas les bords). Nous ne donnons pas une preuve 
compl\`ete de ce th\'eor\`eme, au sens o\`u le proc\'ed\'e de recollement 
utilis\'e est juste esquiss\'e. Nous renvoyons le lecteur au cas du disque 
pour des d\'etails.
\begin{thm}\ \\
Soit $\Gamma$ un arbre fini, connexe et $\mathcal H:=(r_0,d_0,\epsilon,m,h)$ 
une donn\'ee de Hurwitz sur $\Gamma$. Alors $(\Gamma,\mathcal H)$ provient 
d'un automorphisme $\sigma$ d'ordre $p$ d'une 
couronne formelle (ne permutant pas les bords) si et seulement si il 
v\'erifie les conditions $C[i]$, $1 \le i \le 4$ ci-dessous.
\newline\newline
$C[1]$ La racine $r_0$ est de valence $1$
\newline\newline
$C[2]$ Il existe un unique sommet maximal $r_0'$ de $\Gamma_\sigma$ 
tel que l'unique ar\^ete 
$a$ d'origine $r'_0$ v\'erifie $\epsilon(a) \ne 0$. On appelle alors la 
cha\^ \i ne reliant $r_0$ \`a $r'_0$ la cha\^ \i ne fondamentale, et on la 
notera $\Gamma_0$. Les sommets 
(resp. ar\^etes) de la cha\^ \i ne fondamentale seront dits fondamentaux 
(resp. fondamentales).
\newline\newline
$C[3]$ Si $s$ est un sommet maximal de $\Gamma_\sigma$ distinct de $r_0'$, 
l'unique ar\^ete de sommet terminal $s$ est une feuille.
\newline\newline
$C[4]$ Tout sommet de valence sup\'erieure ou \'egale \`a $3$ est 
r\'ealisable.
\newline
\par
De plus, l'\'epaisseur de la couronne formelle est alors la somme 
$\sum_{a \in \aret^+(\Gamma_0)}\epsilon(a)$ des longueurs des ar\^etes 
fondamentales..
\end{thm}
\begin{proof}\ \\
\par
Soit $\sigma$ un automorphisme d'une couronne formelle, ne permutant pas les 
bords. 
$C[1]$ r\'esulte de la construction du mod\`ele minimal. Les assertions 
restantes ont d\'ej\`a \'et\'e vues.
\par
Montrons maintenant la r\'eciproque. On suppose donc que 
$(\Gamma,\mathcal H)$ v\'erifie les propri\'et\'es $C[i]$, $1 \le i \le 4$. 
On voit alors imm\'ediatement le
\begin{lem}\label{lemmededecoupage}
Si $s$ est un sommet fondamental et $a$ une ar\^ete non fondamentale 
d'origine $s$, l'arbre $(\Gamma[a],\mathcal H[a])$ v\'erifie les 
propri\'et\'es $D[i]$, $1 \le i \le 3$ du th\'eor\`eme \ref{hurwitzdisc}.
\end{lem}
On proc\`ede alors par r\'ecurrence sur le nombre $N$ de sommets de la 
cha\^ \i ne fondamentale.
\par
Supposons d'abord $N=2$. L'arbre de Hurwitz est alors r\'eduit \`a sa 
cha\^ \i ne fondamentale, d'apr\`es $C[1]$ et $C[2]$ ($r'_0$ est maximal). 
Notons $a$ l'unique ar\^ete positive de $\Gamma$, et $e:=\epsilon(a)$, 
qui est non nul d'apr\`es $C[2]$. 
Si la racine $r_0$ est additive ou \'etale, alors $m(a)$ est premier \`a $p$. 
Soit $n$ l'entier d\'etermin\'e par $v_K(p)=d_0+n(p-1)$, $Z$ une 
coordonn\'ee de Laurent sur $\mathcal C_e$ ; l'automorphisme $\sigma$ 
de $\mathcal C_e$ donn\'e par $\sigma(Z)=\zeta^{-\frac{1}{m(a)}}Z
(1+\pi^nZ^{m(a)})^{-\frac{1}{m(a)}}$ 
a alors un arbre de Hurwitz \'equivalent \`a $(\Gamma,\mathcal H)$. 
Si maintenant 
$r_0$ est multiplicatif, $h(a) \ne 0 \mod p$. L'automorphisme $\sigma$ 
de $\mathcal C_e$ donn\'e par $\sigma(Z)=\zeta^{-\frac{1}{h(a)}}Z$ 
a alors un arbre de Hurwitz \'equivalent \`a $(\Gamma,\mathcal H)$.
\par
Supposons maintenant $N \ge 3$, et le th\'eor\`eme d\'emontr\'e pour $N-1$. 
Soit $s$ le sommet fondamental qui est l'origine de l'unique ar\^ete $a'$ de 
sommet terminal $r'_0$. Soit 
$\Gamma_1$ le sous-arbre de $\Gamma$ obtenu en retirant les ar\^etes positives 
d'origine $s$, et $\mathcal H_1$ la donn\'ee de Hurwitz restreinte de 
$\mathcal H$ sur $\Gamma_1$. Alors, la cha\^ \i ne fondamentale de 
$(\Gamma_1,\mathcal H_1)$ poss\`ede $N-1$ sommets, et 
$(\Gamma_1,\mathcal H_1)$ v\'erifie les conditions $C[i]$, pour 
$1 \le i \le 4$. Il existe alors un automorphisme $\sigma_1$ d'ordre $p$ de 
la couronne formelle d'\'epaisseur 
$e_1:=\sum_{a \in \aret^+ \Gamma_0, a \ne a'}\epsilon(a)=e-\epsilon(a')$. 
dont l'arbre de Hurwitz est \'equivalent \`a $(\Gamma_1,\mathcal H_1)$.
\par
Par ailleurs, soit $\Gamma_2$ le sous-arbre  de $\Gamma$ de sommets $s$ et 
$r'_0$ et dont l'unique ar\^ete est $a$, et $\mathcal H_2$ la donn\'ee de 
Hurwitz restreinte de 
$\mathcal H$ sur $\Gamma_2$. D'apr\`es le cas $N=2$, il existe 
un automorphisme $\sigma_2$ d'ordre $p$ de 
la couronne formelle d'\'epaisseur 
$e_2=\epsilon(a')$ 
dont l'arbre de Hurwitz est \'equivalent \`a $(\Gamma_2,\mathcal H_2)$.
Par hypoth\`ese, $s$ est r\'ealisable. 
Soit $(u_s,j_s)$ v\'erifiant les conditions de $[M]$ si 
$s$ est multiplicatif, et celles de $[A]$ si $s$ est additif (Le sommet $s$ 
ne peut pas \^etre \'etale). On construit \`a l'aide de $u_s$ et du lemme 
\ref{lemmededecoupage} un automorphisme d'ordre $p$ de $R\{Z,Z^{-1}\}$, et 
on peut alors recoller d'un c\^ot\'e la couronne d'\'epaisseur $e_1$ munie de 
$\sigma_1$ et de l'autre la couronne d'\'epaisseur $e_2$, munie de 
l'automorphisme $\sigma_2$. On obtient alors une couronne formelle 
d'\'epaisseur $e$ munie d'un automorphisme $\sigma$ dont l'arbre de Hurwitz 
est \'equivalent \`a $(\Gamma,\mathcal H)$.
\end{proof}
\subsection{Un crit\`ere suffisant de r\'ealisabilit\'e}

Nous ne donnons ici que le crit\`ere pour un sommet fondamental $s$ 
pour lequel les deux ar\^etes $a_1$ et $a_2$ fondamentales d'origine $s$ 
v\'erifient $m_1:=-m(a_1)>0$ et $m_2:=-m(a_2)>0$. Les autres sommets peuvent 
se traiter par le crit\`ere \ref{critrealdisc}. La preuve est tr\`es similaire 
\`a celle de ce dernier, et est laiss\'ee au lecteur \`a titre d'exercice.
\begin{prop}\ \\
\label{criterekuro}
(a) Supposons que $s$ est multiplicatif. Notons 
$I_s:=\aret(s) \setminus \{a_1,a_2\}$ et 
$\underline h:=(h(a))_{a \in I_s}$. 
Alors, si il existe une 
partition maximale $\underline h$-adapt\'ee $P$ de $I_s$ avec 
$\card(P)\le [\displaystyle{\frac {m_1}p}]+
[\displaystyle{\frac {m_2}p}]+1$, le 
sommet $s$ est r\'ealisable.\\
(b) Supposons que $s$ est additif. 
Notons $I_s$ la r\'eunion disjointe de l'ensemble fini 
$I'_s:=\aret(s) \setminus \{a_1,a_2\}$ 
et d'un ensemble fini \`a 
$m_1+m_2-\card(I'_s)$ \'el\'ements, et 
$\underline e=(e_i)_{i \in I_s}$ la famille d\'efinie par 
$e_i:=m(i)$ si $i \in I'_{s}$ et $e_i=-1$ sinon. 
Alors, si il existe une 
partition maximale $\underline e$-adapt\'ee $P$ de $I_s$ avec 
$\card (P)\le [\displaystyle{\frac {m_1}p}]+
[\displaystyle{\frac {m_2}p}]+1$, le 
sommet $s$ est r\'ealisable.
\end{prop}

\subsection{Structure des automorphismes d'ordre $p$ 
\`a petits conducteurs}
\par
Dans ce qui suit, $\sigma$ d\'esigne un automorphisme d'ordre $p$ de 
$\mathcal C_e$, 
qui ne permute pas les bords, op\'erant sans inertie aux bords. On fera 
l'hypoth\`ese ici que $m_1<p$ et $m_2<p$, o\`u $m_1+1$ et $m_2+1$ d\'esignent 
les conducteurs aux bords $\eta_1$ et $\eta_2$.
\begin{figure}[htpb]
         \input{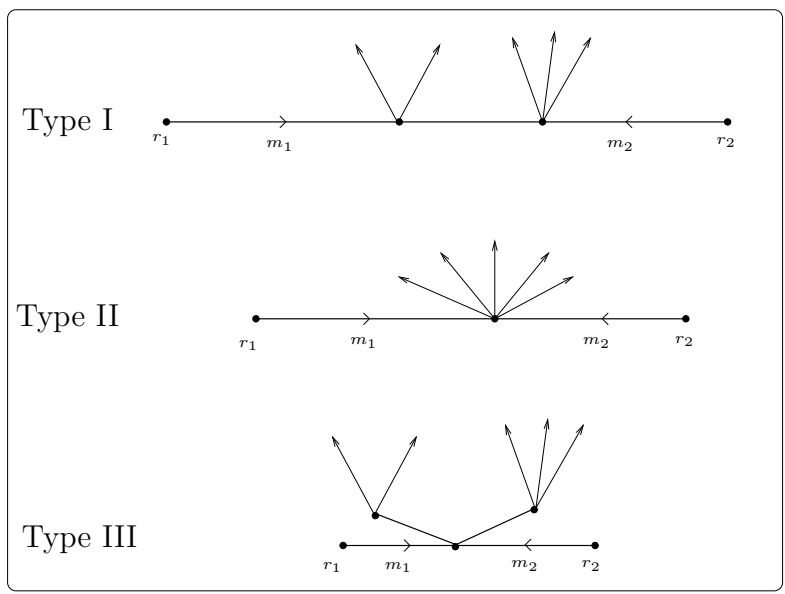}
         \caption{Les trois types d'arbre de Hurwitz de couronnes \`a petits conducteurs}
\end{figure}

\begin{thm}
\label{thmpetitscond}
Soit $\sigma$ un automorphisme de $\mathcal C_e$ d'ordre $p$, ne permutant 
pas les bords, qui op\`ere sans inertie aux bords. Soit
$m_i$ le conducteurs au bord $\eta_i$, $i=1,2$. On note $r_i$ le sommet 
de l'arbre de Hurwitz $(\Gamma_{\sigma}, \mathcal H_{\sigma})$ qui 
correspond au bord $\eta_i$. 
Si $m_1<p$ et $m_2<p$, la cha\^ \i ne fondamentale de 
$(\Gamma_{\sigma}, \mathcal H_{\sigma})$ contient 
au plus deux sommets en dehors de $r_1$, $r_2$. De plus, s'il y en a deux, 
ils sont multiplicatifs. Autrement dit, l'arbre de Hurwitz 
$(\Gamma_\sigma,\mathcal H_\sigma)$ (avec pour fixer les id\'ees la racine 
$r_1$)
est de l'un des trois types ci-dessous : 
\newline
\item 
\begin{center}\bf Type I : \rm $e>(\frac 1{m_1}+\frac 1{m_2})v_K(\zeta -1)$
\end{center}
\par
Si $e>(\frac 1{m_1}+\frac 1{m_2})v_K(\lambda)$, la cha\^ \i ne 
fondamentale de l'arbre de Hurwitz 
$(\Gamma_\sigma,\mathcal H_\sigma)$ est form\'ee de la racine $r_1$, origine 
d'une ar\^ete $a_1$ de sommet terminal multiplicatif $s_1$, lui-m\^eme 
origine d'une ar\^ete dont le sommet terminal $s_2$ est multiplicatif et 
enfin du sommet $r_2$ (sommet terminal de $a_2$). De plus, pour $i=1,2$, 
il y a exactement $m_i$ feuilles 
d'origine $s_i$. Enfin, 
$\epsilon(a_i)=\displaystyle{\frac 1{m_i}}v_K(\lambda)$.
%\begin{center}\includegraphics[width=9cm]{typeI.eps}
%\end{center}
\item
\begin{center}\bf Type II : \rm $e=(\frac 1{m_1}+\frac 1{m_2})v_K(\zeta -1)$
\end{center}
\par
Si $e=(\frac 1{m_1}+\frac 1{m_2})v_K(\lambda)$, la cha\^ \i ne 
fondamentale de l'arbre de Hurwitz 
$(\Gamma_\sigma,\mathcal H_\sigma)$ est form\'ee d'un sommet multiplicatif 
$s$, reli\'e \`a $r_i$ par une ar\^ete $a_i$ au sommet $r_i$, pour  $i=1,2$. 
De plus $s$ est le sommet d'exactement $m_1+m_2$ feuilles. Enfin, 
$\epsilon(a_i)=\displaystyle{\frac 1{m_i}}v_K(\lambda)$.

\item
\begin{center}\bf Type III : \rm $e<(\frac 1{m_1}+\frac 1{m_2})v_K(\zeta -1)$
\end{center}
\par
Si $e<(\frac 1{m_1}+\frac 1{m_2})v_K(\lambda)$, la cha\^ \i ne 
fondamentale de l'arbre de Hurwitz 
$(\Gamma_\sigma,\mathcal H_\sigma)$ est form\'ee d'un sommet additif 
$s$, reli\'e \`a $r_i$ par une ar\^ete $a_i$ au sommet $r_i$, pour  $i=1,2$. 
Le nombre d'ar\^etes non fondamentales d'origine $s$ est inf\'erieur \`a 
$\inf(m_1,m_2)$. Si $a$ est une telle ar\^ete, et si $m(a)<p$, les 
ar\^etes positives d'origine le sommet 
terminal de $a$ sont des feuilles. Il y a au plus une ar\^ete $a$ 
d'origine $s$ pour laquelle $m(a)>p$. Enfin, 
$\epsilon(a_i)<\displaystyle{\frac 1{m_i}}v_K(\lambda)$ pour $i=1,2$.

\end{thm}

\begin{proof}\ \\
\par
L'automorphisme $\sigma$ poss\`ede $m_1+m_2$ points fixes. Ainsi, la 
cha\^ \i ne fondamentale est form\'ee d'au moins trois sommets. 
Soit $s$ un sommet fondamental diff\'erent de $r_1$ et $r_2$, on note $a_1$ 
l'unique ar\^ete fondamentale positive d'origine le sommet terminal $s$ et 
$a_2$ l'unique ar\^ete fondamentale positive d'origine $s$. 
Supposons $m(a_1)>0$ et $m(a_2)>0$. En particulier, $s$ est un sommet additif 
r\'ealisable, on prend les notations de la d\'efinition 
\ref{dfnrealisable}, $[A]$ . Soit 
$\omega=d\bar u_s$. On peut supposer $j_s(\bar a_1)=+\infty$ et $t_{a_2}=0$.  
Comme $\ord_{\infty}d\bar u_s=m(a_1)>0$, quitte \`a ajouter une puissance 
$p$-i\`eme, on peut de plus supposer $\bar u_s$ de la forme 
$$\bar u_s=\frac{P(t)}{\prod_{a \in \aret^+(s)}(t-t_a)^{m(a)}}$$
avec $\deg(P)=\sum_{a \in \aret^+(s)}m(a)$, et $P(t_a) \ne 0$ pour $a$ dans 
$\aret^+(s)$. 
Prenant la coordonn\'ee $x:=t^{-1}$, en notant $Q$ le polyn\^ome 
$Q(x):=x^{\sum_{a \in \aret^+(s)}m(a)}P(\frac 1x)$, on a alors 
$$\bar u_s=\frac{Q(x)}{\prod_{a \in \aret^+(s)}(1-t_ax)^{m(a)}}.$$
Comme par ailleurs $1\le 1+\ord_\infty\omega=m(a_1)\le m_1<p$, 
on a le d\'eveloppement en s\'erie de Laurent 
$$\bar u_s=1+\alpha x^{m(a_1)} \mod x^{m(a_1)+1},$$ 
o\`u $\alpha$ est non nul. Mais alors,
$$Q(x)=(1+\alpha x^{m(a_1)}+x^{m(a_1)+1}l(x))
\prod_{a \in \aret^+(s)}(1-t_ax)^{m(a)}.$$
On a $(\sum_{a \in \aret^+(s)}m(a))-(m(a_1))\le
\sum_{a \in \aret^+(s)}m(a)-m(a_1)=2-\card(\aret(s))<0$. 
Le coefficient de $x^{m(a_1)}$ dans $Q(x)$ est donc $\alpha \ne 0$. Or, 
le degr\'e de $Q$ est \'egal \`a $\sum_{a \in \aret^+(s)}m(a)$, 
on obtient donc une contradiction. On ne peut pas 
avoir $m(a_1)>0$ et $m(a_2)>0$. Il existe donc au plus une ar\^ete 
fondamentale positive $a$ telle que $m(a)>0$, et c'est donc forc\'ement 
l'ar\^ete d'origine $r_1$. De m\^eme, la seule ar\^ete fondamentale positive 
$a$ telle que $m(a)<0$ est l'ar\^ete de sommet terminal $r_2$. On voit qu'il y 
a trois cas : 
Notons $s_1$ le sommet terminal de l'ar\^ete $a_1$ d'origine $r_1$, et 
$a_0$ l'ar\^ete fondamentale positive d'origine $s_1$.
\par
Si $m(a_0)=0$, $s_1$ est multiplicatif, et le sommet 
terminal $s_2$ de $a$ est multiplicatif. On tombe sur le type I. On a 
$m_1=m(a_1)$, les ar\^etes 
positives d'origine $s$ distinctes de $a$ sont des feuilles, et la relation 
$\sum_{a \in \aret(s)}(m(a)+1)=2$ montre qu'il y en a $m_1$. De m\^eme, 
il y a $m_2$ feuilles d'origine $s_2$ (le sommet terminal de $a_0$). Le calcul 
de $\epsilon(a_i)$ se fait en utilisant la loi de variation de la 
diff\'erente. On a alors $e>\epsilon(a_1)+\epsilon(a_2)=(\frac 1{m_1}+
\frac 1{m_2})v_K(\lambda)$.
\par 
Supposons \`a pr\'esent $m(a_0)<0$ et $s_1$ multiplicatif. Alors le sommet 
terminal de $a_0$ est $r_2$. On tombe sur le type II. On voit alors que les 
sommets non fondamentaux d'origine $s_1$ sont des feuilles, et la relation 
$\sum_{a \in \aret(s)}(m(a)+1)=2$ montre qu'il y en a $m_1+m_2$. On voit 
alors facilement que $e=\epsilon(a_1)+\epsilon(a_2)=(\frac 1{m_1}+
\frac 1{m_2})v_K(\lambda)$.
\par
Supposons enfin $m(a_0)<0$ (i.e. $m(a_0)=-m_2$) 
et $s_1$ additif. Soit $q$ le nombre d'ar\^etes 
non fondamentales d'origine $s=s_1$. On peut alors supposer, avec les 
notations de la proposition \ref{criterekuro},
$$\bar u_s=\frac{P(t)}{\prod_{a \in I'_s}(t-t_a)^{m(a)}}$$
o\`u $\deg(P)=\sum_{a \in I'_s}m(a)$, $j_s(\bar a_1)=\infty$, $j_s(a_2)=0$, 
$j_s(a)=t_a \ne 0$ pour $a$ dans $I'_s$, $P(0)\prod_{a \in I'_s}P(t_a)\ne0$, 
et $\ord_{\infty}d\bar u_s=m_1-1$, 
$\ord_{0}d\bar u_s=m_2-1$. On a $\bar u_s=1+\alpha t^{m_2} \mod t^{m_2+1}$, 
o\`u $\alpha \ne 0$. Ainsi, 
$$P(t)=(1+\alpha t^{m_2}+t^{m_2+1}l(x))
\prod_{a \in I'_s}(1-t_a^{-1}t)^{m(a)}.$$
Or, $\deg(P)=\sum_{a \in I'_s}m(a)=m_1+m_2-q$. Si $m_1<q$, on a $\deg(P)<m_2$ 
et on obtient alors une contradiction en remarquant que le coefficient de 
$t^{m_2}$ doit \^etre $\alpha$. Donc $q \le m_1$, et en regardant \'egalement 
le d\'evelloppement de $\bar u_s$ en $\infty$, on voit que $q \le m_2$. 
Le fait qu'il y ait au plus une ar\^ete non 
fondamentale $a$ d'origine $s$ pour laquelle $m(a)>p$ 
r\'esulte de $m_1+m_2<2p$. Si $a$ est une ar\^ete non 
fondamentale d'origine $s$ pour laquelle $m(a)<p$, toute ar\^ete 
d'origine $t_{\Gamma_\sigma}(a)$ est une feuille 
d'apr\`es le th\'eor\`eme III 3.1 de \cite{G-M 2}.
 Enfin, la loi de 
variation de la diff\'erente entra\^ \i ne que $e<(\frac 1{m_1}+
\frac 1{m_2})v_K(\lambda)$.

%\begin{center}\includegraphics[width=14cm]{types.eps}
%\end{center}
\end{proof}

\end{document}